\documentclass{article}
\usepackage{stmaryrd}
\usepackage{amsfonts}
\usepackage{amsmath}
\usepackage[all]{xy}
\newtheorem{theorem}{Theorem}[section]
\newtheorem{definition}[theorem]{Definition}
\newtheorem{proposition}[theorem]{Proposition}

\newtheorem{lemma}[theorem]{Lemma}

\usepackage{mathrsfs}
\usepackage{amssymb}
\usepackage{amscd}

\date{}

\usepackage{amssymb}
\usepackage{extarrows}
\usepackage{cite}

\begin{document}

\title{Construction of   free   Lie Rota-Baxter superalgebra via  Gr\"{o}bner-Shirshov bases theory}

\author{
 Jianjun Qiu  \\
{\small \ School of Mathematics and Statistics, Lingnan Normal
University}\\
{\small Zhanjiang 524048, P. R. China}\\
{\small jianjunqiu@126.com}\\
\\
Yuqun Chen\footnote{Corresponding author}  \\
{\small \ School of Mathematical Sciences, South China Normal
University}\\
{\small Guangzhou 510631, P. R. China}\\
{\small yqchen@scnu.edu.cn}\\
}

\maketitle \noindent\textbf{Abstract:} In this paper,  we construct  free   Lie Rota-Baxter superalgebra by using   Gr\"{o}bner-Shirshov bases theory.  We firstly construct free  operated  Lie  superalgebras by the  operated  super-Lyndon-Shirshov monomials.  Secondly,    we establish    Gr\"{o}bner-Shirshov bases theory for  operated   Lie  superalgebras.    Thirdly,   we  find a Gr\"{o}bner-Shirshov basis of  a   free  Lie  Rota-Baxter superalgebra on a $\mathbb{Z}_2$-graded set.  Consequently,  we can obtain      a linear basis of a  free  Lie  Rota-Baxter superalgebra by the composition-diamond lemma for   operated   Lie  superalgebras.

 \ \


  \ \


\tableofcontents

\section{Introduction}

Gr\"{o}bner bases and Gr\"{o}bner-Shirshov bases  were invented
independently by  Shirshov for ideals of  free
Lie algebras \cite{Sh,Shir3} and implicitly free associative
algebras \cite{Sh,Shir3}  (see also \cite{Be78,Bo76}), by
Hironaka \cite{Hi64} for ideals of the power series algebras (both
formal and convergent), and by  Buchberger \cite{Bu70} for ideals
of the polynomial algebras. Gr\"{o}bner bases and
Gr\"{o}bner-Shirshov bases theories have been proved to be very
useful in different branches of mathematics, including commutative
algebra and combinatorial algebra. It is a powerful tool to solve
the following classical problems: normal form; word problem;
conjugacy problem; rewriting system; automaton; embedding theorem;
PBW theorem;  extension; homology;  growth function; Dehn function;
complexity; etc. See, for example, the books \cite{AL, BKu94, BuCL,
BuW, CLO, Ei}, the papers \cite{Be78,   Bo76,       Mikhalev92}, and the
surveys \cite{BC14,     BK03}.

  Let $(R, \cdot)$ be an algebra over a field  $k$ and  $\lambda\in k$.  A  linear operator   $ P:R\rightarrow R$   satisfying
$$
 P(x)\cdot P(y)  = P(  P(x)\cdot y)  +P( x\cdot P(y))+\lambda P( x\cdot y ),  \   \forall\  x,y \in  R
$$
 is called a   Rota-Baxter   operator on the algebra $(R, \cdot)$  of weight $\lambda$.  An algebra $(R, \cdot)$  together with   a $\lambda$-Rota-Baxter   operator $P$ on it   is called  a   Rota-Baxter   algebra of weight $\lambda$ \ (or $\lambda$-Rota-Baxter   algebra).

The Rota-Baxter operators on an associative algebra  were introduced  by    Baxter \cite{Bax60}   to solve an analytic  problem in probability  and    was studied by    Rota  \cite{Ro69} in combinatorics later. The research results  on  Rota-Baxter operators  on an associative algebra as well as   associative  Rota-Baxter algebras are  very rich (see \cite{lguo} and the references given there).    As we know, the  free object of a variety   of   algebras  play a great   role in the theories  and applications.    There are many constructions of free  associative $\lambda$-Rota-Baxter algebras on both commutative and noncommutative cases   by using  different methods, for example, see \cite{am2006,  bcq2010, Ca72,    EG08a, EG08b,    lguo,    GK00a, GK00b,  Ro69}.

The  Rota-Baxter operators   on a Lie algebra  are also   called the operator form of the classical Yang-Baxter equation due to Semenov-Tian-Shansky's work \cite{sts}.
There are  some new results on  Rota-Baxter operators on Lie algebras, Lie  Rota-Baxter algebras  and related topic, for example,   see \cite{amm2015, ab, bgn10, lhb, pgb, gk2016}.  Recently,    free Lie Rota-Baxter algebras were constructed by Gubarev \cite{gub2016},  V. Gubarev and     Kolesnikov \cite{gk2016},     Qiu and Chen \cite{qc2017}.

The concept of   Rota-Baxter operators   on a Lie superalgebra   as well as  Lie Rota-Baxter superalgebras  were appeared in the recent  papers \cite{amm2015,  whc2010}, which are  also closed to the Yang-Baxter equations in the Lie superalgebras case.
  A Lie superalgebra $L=L_{\bar{0}}\oplus L_{\bar{1}}$  equipped with an  even linear map $ P:L\rightarrow L$ satisfying
$$
[P(x)P(y)]  = P([P(x) y] )-(-1)^{|x||y|}P([P(y)x])+\lambda P ( [x y ])
$$
for all homogenous elements $x$ and $y$,   is called a  Lie   Rota-Baxter  superalgebra of weight $\lambda$ (or Lie $\lambda$-Rota-Baxter  superalgebra).

  Rota-Baxter algebras are algebras with linear operators satisfying some identities. An algebra with operators  are called operated algebras  \cite{lguo2009} or  $\Omega$-algebras \cite{Ku60}.
In  \cite{bcq2010}, Bokut, Chen and Qiu established   Gr\"{o}bner-Shirshov bases theory for associative  algebras with operators, which was used to give a linear basis of a  free associative Rota-Baxter algebra on a set. Recently,   Qiu and  Chen \cite{qc2017} applied    Gr\"{o}bner-Shirshov bases theory for  Lie   algebras  with  operators   to construct  free   Lie  Rota-Baxter algebras.  In this paper, based on  Gr\"{o}bner-Shirshov bases theory for  Lie  superalgebras developed   by  Mikhalev \cite{Mikhalev92},
   we establish    Gr\"{o}bner-Shirshov bases theory for  Lie  superalgebras with operators, which   is   applied     to construct   free  Lie Rota-Baxter superalgebras.   Especially, we  give a linear basis of  a free    Lie Rota-Baxter superalgebra  on a  $\mathbb{Z}_2$-graded set.

The rest of the  paper is organized as follows. In Section 2, we  review the properties of  super-Lyndon-Shirshov  words and  monomials.  In Section 3, we  develop the  Gr\"{o}bner-Shirshov
bases   for   operated  associative     superalgebras.  In Section 4, we construct    free  operated  Lie  superalebras    by the   operated   super-Lyndon-Shirshov  monomials.     In Section 5, we  give definition of Gr\"{o}bner-Shirshov
bases   for   operated  Lie   superalgebras and establish Composition-Diamond lemma for operated   Lie superalgebras. In Section 6, we find   a Gr\"{o}bner-Shirshov basis of   a free   Lie  Rota-Baxter superalgebra, which is applied to give    a  linear basis of such  algebra    by Composition-Diamond lemma for  operated  Lie  superalgebras.

In the rest of this paper, unless otherwise specified, $k$ is a   field with characteristic zero. We will denote by  $kY$    the $k$-linear space generated by  the set $Y$.  By a $\mathbb{Z}_2$-graded set we mean a disjoint union  $X=X_{\bar{0}}\cup  X_{\bar{1}}$ of two subsets.

\section{Free  Lie  superalgebras}

\subsection{Basic definitions of  superalgebras}

In this subsection, we review some basic concepts and  properties of   Lie superalgebras.  For more details we refer the reader to \cite{ka1977,sch1979}.

\begin{definition} (\cite{ka1977,sch1979})
(1)\  Let   $V$ be a linear space over  the field  $k$. If $V=V_{\bar{0}}\oplus V_{\bar{1}}$, then $V$ is called  a  $\mathbb{Z}_2$-graded  space.

(2) \ A subspace $U$ of $V$ is called $\mathbb{Z}_2$-graded, if  $U=(U\cap V_{\bar{0}})\oplus (U\cap V_{\bar{1}})$.

(3) \ The elements of  $V_{\bar{0}}$ (resp. $V_{\bar{1}}$) are said to be
even  (resp. odd).    The elements of $V_{\bar{i}}, \ \bar{i} \in \mathbb{Z}_2=\{\bar{0}, \bar{1}\}$  are said to be homogenous with
parity $\bar{i}$.  The parity of a homogeneous element $x$ is denoted by $|x|$.
\end{definition}

\begin{definition}  (\cite{ka1977,sch1979})  (1)\
Let $X=X_{\bar{0}}\cup X_{\bar{1}}$ be a  $\mathbb{Z}_2$-graded and  $Y=Y_{\bar{0}}\oplus Y_{\bar{1}}$  a   $\mathbb{Z}_2$-graded  vector space.    A map $\phi: X\rightarrow Y$ is said to be even (resp. odd) if $\phi(X_{\bar{i}})\subseteq Y_{\bar{i}}$ (resp. $\phi(X_{\bar{i}})\subseteq W_{\bar{i}+\bar{1}})$ for any  $\bar{i}\in \mathbb{Z}_2$.

(2)\ Let $V=V_{\bar{0}}\cup V_{\bar{1}}$ and $W=W_{\bar{0}}\oplus W_{\bar{1}}$  be two   $\mathbb{Z}_2$-graded vector spaces.        A  linear map $\phi: V\rightarrow W$ is said to be even (resp. odd) if $\phi(V_{\bar{i}})\subseteq W_{\bar{i}}$ (resp. $\phi(V_{\bar{i}})\subseteq W_{\bar{i}+\bar{1}})$ for any  $\bar{i}\in \mathbb{Z}_2$.

(3)\ If  a  linear map $\phi: V\rightarrow W$  is  even, then  $\phi$  is also said to be a homomorphism from the  $\mathbb{Z}_2$-graded space $V $ to  $W$.
\end{definition}

\begin{definition} (\cite{ka1977,sch1979})
(1)  \ An algebra  $(A, \circ)$    over the field $k$ is said to be a superalgebra if   $A=A_{\bar{0}}\oplus A_{\bar{1}}$ and $A_{\bar{i}}\circ A_{\bar{j}}\subseteq A_{\bar{i}+\bar{j}}$ for any $\bar{i}, \ \bar{j} \in \mathbb{Z}_2$.

(2)\  A super-subalgebra  (resp. super-ideal) of a superalgebra $A=A_{\bar{0}}\oplus A_{\bar{1}}$ is a subalgebra  (resp. ideal) of $A$ which is, in addition a $\mathbb{Z}_2$-graded subspace of the $\mathbb{Z}_2$-graded    space $A$.

(3)  A homomorphism of superalgebras is a homomorphism of the underlying algebras as well as of the underlying $\mathbb{Z}_2$-graded   space.
\end{definition}

\begin{definition} (\cite{ka1977,sch1979})
(1)\ A superalgebra $(A, \cdot)$ is called an associative superalgebra if for all homogenous elements $x, y, z$,
$$
x\cdot(y\cdot z)=(x\cdot y)\cdot z.
$$

(2) \ A superalgebra $(L, [  ,   ])$  is called a Lie  superalgebra if for all homogenous elements $x, y, z$,
$$
[x,y]=-(-1)^{|x||y|}[y,x],\ \ \  \mbox{(super-skew-symmetry)},
$$
$$
[x,[y,z]]=[[x,y],z]+(-1)^{|x||y|}[y,[x,z]], \ \mbox{(super-Jacobi identity)}.
$$
\end{definition}

Let $(A, \cdot)$ be an associative superalgebra, where $A=A_0\oplus A_1$.  It is well known that  $(A, [, ])$ is a  Lie superalgebra, where $[, ]$ is defined by
$$
[a, b]=a\cdot b-(-1)^{|a||b|}b\cdot a
$$
for all homogenous elements $a$ and $b$.

\subsection{Definitions of super-Lyndon-Shirshov words and monomials}

In this subsection, we review the concept   of Lyndon-Shirshov (LS) and  super-Lyndon-Shirshov (SLS) words and monomials, which can be found in \cite{Bo1999, Sh, Shir3, BMP2011,  Mikhalev92}.

For any set $\mathbb{Z}_2$-graded set $X=X_{\bar{0}}\cup X_{\bar{1}}$, let $S(X)$ (resp. $X^*$) be the  free semigroup (resp. monoid) on $X$ of nonempty   associative words (resp.  associative  words with empty word 1).      For any $w=x_1x_2\cdots x_n\in S(X)$, where each  $x_i\in X$,
define the parity of $w$ by
$
|w|=\sum_{i=1}^n|x_i|,
$
where $|x|=\bar{i}$, if $x\in X_{\bar{i}}, \ \ i=0,1$.
By   the parity of the words in $S(X)$,  $S(X)$  has the $\mathbb{Z}_2$-grading
$$
S(X)=S(X)_{\bar{0}}\cup S(X)_{\bar{1}},
$$
where
$S(X)_{\bar{0}}=\{u\in S(X)| |u|= \bar{0}\}$ and $S(X)_{\bar{1}}=\{u\in S(X)| |u|= \bar{1}\}$.

Let us denote by $X^{**}$ the free groupoid of     nonassociative words  on  $X$.  Assume that $X=X_{\bar{0}}\cup X_{\bar{1}}$ is a  $\mathbb{Z}_2$-graded set.   For any nonassociative word $(u)\in X^{**}$, we define the  parity of $(u)$ by $|(u)|=|u|$. For example, if $w=x_0x_0x_1$ and $(w)=(x_0(x_0x_1))$, where $x_0\in X_{\bar{0}}$ and $x_1\in X_{\bar{1}}$,   then    $|(w)|=|w|=\bar{1}$.
Let  $X^{**}_{\bar{0}}=\{(u)\in X^{**}| |(u)|= \bar{0}\}$ and $X^{**}_{\bar{1}}=\{(u)\in X^{**}| |(u)|= \bar{1}\}$. It follows that  $X^{**}$ is also a    $\mathbb{Z}_2$-graded set, i.e.  $X^{**}=X^{**}_{\bar{0}}\cup X^{**}_{\bar{1}}$.

For any $u\in X^*$, let $deg(u)$ to be  the degree (length) of $u$.  Let $>_{X}$  be  a   well order on $X$.    Define the lex-order $>_{lex}$ and the deg-lex order $>_{dl}$ on $X^*$ with respect to $>_{X}$ by:

 (i) $1>_{lex} u$ for  any nonempty word  $u$, and  if $u=x_i\tilde{u} $ and $v=x_j\tilde{v}$, where $x_i, x_j\in X$,   then $u>_{lex} v$ if   $x_i> x_j$, or $x_i=x_j$ and $\tilde{u}>_{lex} \tilde{v}$ by induction.

 (ii) $u>_{dl}v$ if $deg(u)>deg(v)$, or $deg(u)=deg(v)$ and $u>_{lex}v$.

\begin{definition}(\cite{Bo1999, Sh, Shir3, BMP2011,  Mikhalev92})
A nonempty associative  word $u$ is called  a  Lyndon-Shirshov (LS)  word on $X$ with respect to the   order $>_{lex}$,  if $  u=vw >_{lex} wv$ for any  decomposition of $u=vw$, where $v, w\in S(X)$.
\end{definition}

\begin{definition}(\cite{Bo1999, Sh, Shir3, BMP2011,  Mikhalev92})
A nonassociative word $(u)\in X^{**}$ is said to be   a   Lyndon-Shirshov (LS)  monomial  on $X$ with respect to the    order $>_{lex}$,     if
  \begin{itemize}
    \item [(a)] $u$ is an  LS word on $X$;
    \item [(b)] if $(u)=((v)(w))$, then both $(v)$ and $(w)$ are  LS  monomials on $X$;
    \item [(c)] if $(v)=((v_1)(v_2))$, then $v_2 \leq_{lex} w$.
  \end{itemize}
\end{definition}

\begin{definition}(\cite{Bo1999, Sh, Shir3, BMP2011,  Mikhalev92})
A nonempty associative  word $u$ is called  a    super-Lyndon-Shirshov (SLS)  word on $X$, if it is an  LS word   or it has the form $u=vv$ with $v$ is an LS word and $|v|=\bar{1}$.
\end{definition}

\begin{definition}(\cite{Bo1999,   Mikhalev92})
A   nonassociative word $(u)$ is called an  super-Lyndon-Shirshov (SLS) monomial if either it is a LS monomial or it has the form $(u)=((v)(v))$ with $(v)$ is a  LS monomial and $|(v)|=\bar{1}$.
\end{definition}

\subsection{Properties of super-Lyndon-Shirshov words and monomials}

The following lemma asserts that there is a natural 1-1 correspondence
between the set of SLS words and the set of SLS monomials.

\begin{proposition}\label{pro2.9}  {\em (\cite{Bo1999, BMP2011,  Mikhalev92})}
If $(u)$ is a SLS  monomial,
then $u$ is a SLS word. Conversely, for any SLS word $u$, there is a unique arrangement of brackets $[u]$ on $u$
such that $[u]$  is a SLS monomial.
\end{proposition}

For any LS word $u$, there are some bracketing ways    to  get the unique LS monomial $[u]$.  One is up to down bracketing way, which is defined inductively by
$[x_i]=x_i$ for $x\in X$ and $[u]=[[v][w]]$, where $u=vw$ and $w$ is the longest LS proper subword  in  $u$ from the right hand side.
Moreover, if $u=vv$ is a SLS word with  $v$ is  an LS word and   $|v|=\bar{1}$, then  $[u]=[[v][v]]$.

We will denote by $SLSW(X)$    (resp. $SLSM(X))$  the set of all  SLS words  (resp. monomials) on $X=X_{\bar{0}}\cup X_{\bar{1}}$ with respect to the lex-order $>_{lex}$.  According to  Proposition  \ref{pro2.9}, we can get
$$
SLSM(X)=\{[u] | u\in SLSW(X)\}.
$$

\begin{proposition} {\em (\cite{Bo1999, BMP2011,  Mikhalev92})} \label{pro2.10}  Let $w_1, w_2$  be  SLS words and $w_1=e_1e_2, w_2=e_2e_3$  where $e_1, e_2, e_3$ are nonempty word. Then we have the following three results.

(1)   If $w_1=u$,  $w_2=v^n\ (n=1, 2)$,   $u, v$ are LS word and  $u\neq v$,  then $e_1e_2e_3$ is an LS word.

(2)   If $w_1=uu$,    $w_2=v^n \ (n=1, 2)$,      $u, v$ are LS word and  $u\neq v$, then $e_1=u\tilde{e}_1$ and $e_1e_2e_3$ is an LS word, where $\tilde{e}_1$ may be empty.

(3)   If $w_1=uu$,  $w_2=uu$, and $u$ is an  LS word with $|u|=1$, then $e_1=e_2=e_3=u$.
\end{proposition}

\begin{proposition} {\em (\cite{Bo1999, BMP2011,  Mikhalev92})}   \label{pro2.11} Let $u, v$ be LS words and $u\neq v$.
 If   $v^l$ is a subword of $u^m$, then
 $
 v\neq u_2u^{l_1}a, l_1\geq 0, u_1u_2=u=ab, u_2, a\in S(X)
 $
 and
 $
 u\neq dv^{k_1}v_1, k_1\geq 0, v_1v_2=v=cd, v_1, d\in S(X).
 $
It follows that  $v^l$ is a subword of $u$.
\end{proposition}

\subsection{Free  Lie superalgebras}

\begin{definition} A free   associative (resp. Lie)  superalgebra on a $\mathbb{Z}_2$-graded $X=X_{\bar{0}}\cup X_{\bar{1}}$  is  a pair $(F, \textbf{i})$  with $F$ an  associative (resp. Lie)  superalgebra and $\textbf{i} : X \rightarrow F$ an even map from $X$ to $F$ such that if $L$ is
any (resp. Lie)    superalgebra and $\textbf{j} : X \rightarrow L$ is an even map, then there is a unique homomorphism    of  associative (resp. Lie)   superalgebras  $\phi_\textbf{j}:  F\rightarrow L$ such that  $\phi_\textbf{j}\textbf{i}=\textbf{j}$.
\end{definition}

\begin{theorem}(\cite{Bo1999})\     Let  $\textbf{AS}(X)=\textbf{AS}(X)_{\bar{0}}\oplus  \textbf{AS}(X)_{\bar{1}}$, where
$
\textbf{AS}(X)_{\bar{i}}=kS(X)_{\bar{i}}, i=1,2.
$
Then $\textbf{AS}(X)=\textbf{AS}(X)_{\bar{0}}\oplus  \textbf{AS}(X)_{\bar{1}}$   is a free nounitary associative superalgebra on the   $\mathbb{Z}_2$-graded set $X=X_{\bar{0}}\cup X_{\bar{1}}$.
\end{theorem}

The algebra $(\textbf{AS}(X), [,])$ becomes a  Lie superalgebra with the superbracket
$$
[u, v]=uv-(-1)^{|u||v|}vu,
$$
where  $u$ and $v$ are homogenous elements.  Let $\textbf{LS}(X)$ be the Lie  sub-superalgebra of the Lie superalgebra $(\textbf{AS}(X), [,])$ generated by $X$  under the superbracket $[,]$.

Given a nonzero element $f\in \textbf{AS}(X)$,  we denote by $\bar{f}$ the leading  word
appearing in $f$ with respect to   the  deg-lex order $>_{dl}$ on $X^*$ and $lc(f)$ the    leading coefficient of  $f$.

\begin{lemma}\label{le2.2}{\em (\cite{Bo1999, BMP2011,  Mikhalev92})}
For any  $(u)\in X^{**}$,   $(u)$ has a representation
$$
(u)=\sum \alpha_i[u_i],
$$
in $\textbf{LS}(X)$, where each $\alpha_i\in k$, $[u_i]\in SLSM(X)$ and $deg(u)=deg(u_i)$.
\end{lemma}

\begin{lemma} \label{le2.3}  {\em (\cite{Bo1999, BMP2011,  Mikhalev92})} \label{lem2.2}
For any  SLS word  $u\in SLSW(X)$,
$$
[u]=lc([u])u+\sum \alpha_i[u_i],
$$
where $lc([u])$ is the leading coefficient of $[u]$ in $\textbf{LS}(X)\subseteq \textbf{AS}(X)$,  $u_i<_{Dl} u$  and $deg(u)=deg(u_i)$. If follows that   $\overline{[u]}=u$.
\end{lemma}

\begin{theorem} \label{thm2.4} {\em (\cite{Bo1999, BMP2011,  Mikhalev92})} \label{th2.3}   $\textbf{LS}(X)$ is a  free  Lie superalgebra generated by the $\mathbb{Z}_2$-graded set  $X=X_{\bar{0}}\cup X_{\bar{1}}$. Moreover,
the set $SLSM(X)$ consisting  of  all   SLS monomials  forms a linear basis  of  $\textbf{LS}(X)$.
\end{theorem}

\begin{theorem}\label{th2.15}  {\em (\cite{Bo1999, BMP2011,  Mikhalev92})}  Assume that  $u, v$  are  LS words on $X$.
 If   $u=avb$, then    $[u]=[a[vc]d]$, where $b=cd$ and $c$ may be empty.  Define  an   arrangement of brackets  $[u]_{v}$ on $u$ by
$$
[u]_{v}=[u]_{[vc]\mapsto [[\cdots[[[v][c_1]][c_2]]\cdots ][c_m]]},
$$
where $c=c_1c_2\cdots c_m$ with each  $c_i\in ALSW(X)$, $c_{1} \leq_{{  lex}} c_{2}\leq_{{  lex}} \cdots \leq_{{  lex}} c_m$. Then we have
$$
[u]_v=a[v]b+\sum \alpha_ia_i[v]b_i,
$$
where each $\alpha_i\in k, a_i, b_i\in X^*$ and $a_ivb_i<_{{dl}} avb=u$. It follows that   $\overline{[u]_{v}}=u$.
\end{theorem}

Based on the work  of  Mikhalev in  \cite{Mikhalev92}, we define an arrangement of brackets $[u]'_{vv}$ on  LS  words  of the form  $u=avvb$   with   $v$ is an LS word and $|v|=1$   by induction on the degree of $u$.  Let $x$ be the least letter from $X$ which occurs in $u$.

 If $deg(u)=3$,  then $u=yvv, y\in X, y>v=x$ or $u=vvx, v\in X, v>x$. Define $[u]'_{vv}=[y[vv]]$ if $u=yvv$ and $[u]'_{vv}=[[vv]x]$ if $u=vvx$.

 Suppose that we   have been  defined $[u]'_{vv}$ for $u=avvb$ with $3\leq deg(u)<n$.  Let    $deg(u)=n$.

  Case 1. If $v=x\in X$,     the least letter   occurring  in $u$,  then we have $u=\tilde{a}x_ix^svv\tilde{b}$, where $x_i>x, s\geq 0$.     Note that $x_ix^svv=x_ix^{s+2}$ is an LS word and $[u]=[\tilde{a}[x_ix^svv]\tilde{b}]$. By Theorem \ref{th2.15}, we  can   define
$$
[u]'_{vv}=[u]_{x_ix^{s+2}}|_{[x_ix^{s+2}]\mapsto [[x_ix^s][vv]]}.
$$

 Case 2.  If $v\neq x$, then we may  set   $b=x^m\tilde{b}$, where $m\geq 0$ and the first letter of $\tilde{b}$ is not $x$.

 Case 2-1.   If  $m>0$, then $vvx^m$ is an LS word.   By Theorem \ref{th2.15}, we   can  define
$$
[u]'_{vv}=[u]_{vvx^m}|_{[vvx^m]\mapsto [vv](\textbf{ad}x)^m},
$$
where  $[vv](\textbf{ad}x)^m=[[\cdots[[vv]x]\cdots]x]$.

Case 2-2.  If $m=0$, then the first letters of $u$ and $v$ are differential from $x$. Writing down $u$ and $v$ in the new letters $x_ix^l$, i.e. $u'=x_{i_1}x^{l_1}x_{i_2}x^{l_2}\cdots x_{i_p}x^{l_p}$ and   $v'=x_{i_s}x^{l_s}x_{i_{s+1}}x^{l_{s+1}}\cdots x_{i_t}x^{l_t}$. Thus, we have $u'= a'v'v'b'$.    Consider the following order relation on the new letters: $x_ix^l>x_jx^z$ for $x_i>x_j$,  or $x_i=x_j$ with  $l<z$.  Note that  $u'$ and $v'$ are also  LS words with  $|v'|=1$  in  the new letters $x_ix^l$.  Moreover,   the degree of $u'$  in the new letters $x_ix^l$ is less that  the degree of $u$ in $X$.  By the inductive hypothesis, we have  an arrangement of brackets $[u']'_{v'v'}$   on   $u'$ in   the new letters $x_ix^l$.  Define $[u]'_{vv}$ by replacing  the new letters $x_ix^l$ in  $[u']'_{v'v'}$ by  $[x_ix^l]$, i.e.
$$
[u]'_{vv}=[u']'_{v'v'}|_{x_ix^l\mapsto [x_ix^l]}.
$$

\begin{lemma}\label{lem2.18} Let $v$ and $vvx^m$ be  LS words with $|v|=1, m>0$. Then
$$
[vv](\textbf{ad}x)^m=[vv]x^m+\sum_{\substack{i+j=m\\   0< i \leq m, 0\leq j < m}} \beta_{ij}x^i[vv]x^j,
$$
where  $[vv](\textbf{ad}x)^m=[[\cdots[[vv]x]\cdots]x]$.
\end{lemma}
{\bf Proof.} The proof is by induction on $m$.  If $m=1$, then
$$
[vv]\textbf{ad}x=[[vv], x]=[vv]x-(-1)^{|[vv]||x|}x[vv]=[vv]x-x[vv].
$$
Assume that  the  result is true for $m-1$, i.e.
$$
[vv](\textbf{ad}x)^{m-1}=[vv]x^{m-1}+\sum_{\substack{s+t=m-1\\   0< s \leq m-1, 0\leq t < m-1}} \tilde{\beta}_{st}x^s[vv]x^t,
$$
where each $\tilde{\beta}_{st}\in k$. We will prove it for $m$.  Note that
$$
[vv](\textbf{ad}x)^m=([vv](\textbf{ad}x)^{m-1})\textbf{ad}x=[[vv](\textbf{ad}x)^{m-1}, x].
$$
It follows that
\begin{eqnarray*}
[vv](\textbf{ad}x)^m&=& [[vv]x^{m-1}+\sum_{\substack{s+t=m-1\\   0< s \leq m-1, 0\leq t < m-1}} \tilde{\beta}_{st}x^s[vv]x^t, x]\\
&=& [[vv]x^{m-1},x]+  \sum_{\substack{s+t=m-1\\   0< s \leq m-1, 0\leq t < m-1}}\tilde{\beta}_{st}[x^s[vv]x^t, x]\\
&=& [vv]x^{m}-(-1)^{|x|^m}x[vv]x^{m-1}\\
&& +  \sum_{\substack{s+t=m-1\\   0< s \leq m-1, 0\leq t < m-1}}\tilde{\beta}_{st}( x^s[vv]x^{t+1} -(-1)^{|x|^m} x^{s+1}[vv]x^{t})\\
&=& [vv]x^m+\sum_{\substack{i+j=m\\   0< i \leq m, 0\leq j < m}} \beta_{ij}x^i[vv]x^j,
\end{eqnarray*}
where $\beta_{1,m-1}=-(-1)^{|x|^m}$ and $\beta_{ij}=\tilde{\beta}_{i,j-1}-\tilde{\beta}_{i-1,j}(-1)^{|x|^m}$.  Therefore, the result is true for any $m\geq 1$. \hfill $\square$

\begin{theorem} \label{th2.17}  Assume that  $u=avvb$ and $v$  are  LS words with  $|v|=1$. Then
$$
[u]'_{vv}=\alpha a[vv]b+\sum \alpha_ia_i[vv]b_i,
$$
where each $0\neq \alpha, \alpha_i\in k, a_i, b_i\in X^*$ and $a_ivvb_i<_{{dl}} avvb=u$.  It follows that   $\overline{[u]'_{vv}}=u$.
\end{theorem}
{\bf Proof.} We will prove the result by induction on the degree of $u$. Let $x$ be the least letter from $X$ which occurs in $u$.  If $deg(u)=3$, then    $[u]'_{vv}=[y[vv]]$ if $u=yvv$ and $[u]'_{vv}=[[vv]x]$ if $u=vvx$.  It is easy to check that
 $$
 [u]'_{vv}=[y[vv]]=y[vv]-[vv]y, \ vvy<_{{dl}} yvv=u,
 $$
 and
 $$
 [u]'_{vv}=[[vv]x]=[vv]x-x[vv],  xvv<_{{dl}} vvx=u.
 $$

 Suppose the result is valid for   $[u]'_{vv}$  with  $3\leq deg(u)<n$.
Let    $deg(u)=n$.  There are two cases to consider.

Case 1. If $v=x\in X$,     then we have $u=\tilde{a}x_ix^mvv\tilde{b}=\tilde{a}x_ix^{m+2}\tilde{b}$, where $x_i>x, m\geq 0$.   Moreover, by the definition of  $[u]'_{vv}$, we have
$$
[u]'_{vv}=[u]_{x_ix^{m+2}}|_{[x_ix^{m+2}]\mapsto [[x_ix^m][vv]]}.
$$
According to  Theorem \ref{th2.15},
$$
[u]_{x_ix^mvv}= \tilde{a}[x_ix^mvv]\tilde{b} + \sum_{1\leq l\leq z} \alpha_l\tilde{a}_l[x_ix^mvv]\tilde{b}_l.
$$
A simple computation gives
$
[x_ix^mvv]=\frac{1}{2}[[x_ix^m][vv]],
$
which shows  that
$$
[u]'_{vv}=\frac{1}{2} \tilde{a}[[x_ix^m][vv]]\tilde{b} + \sum_{1\leq l\leq z} \frac{\alpha_l}{2}\tilde{a}_l[[x_ix^m][vv]]\tilde{b}_l.
$$
Similar to Lemma \ref{lem2.18},
$$
[x_ix^m]=x_ix^m+\sum_{\substack{s+t=m\\   0< s \leq m, 0\leq t < m}} \beta_{st} x^sx_ix^t.
$$
which gives that
\begin{eqnarray*}
&& [[x_ix^m][vv]]\\
&=&[x_ix^m+\sum_{\substack{s+t=m\\   0< s \leq m, 0\leq t < m}} \beta_{st} x^sx_ix^t, [vv]]\\
&=&[x_ix^m, [vv]]  +\sum_{\substack{s+t=m\\   0< s \leq m, 0\leq t < m}} \beta_{st} [x^sx_ix^t, [vv]]\\
&=& x_ix^m[vv]-[vv]x_ix^m+\sum_{\substack{s+t=m\\   0< s \leq m, 0\leq t < m}} \beta_{st} (x^sx_ix^t [vv]-[vv]x^sx_ix^t).
\end{eqnarray*}
It follows that
\begin{eqnarray*}
&&[u]'_{vv}\\
&=&\frac{1}{2} \tilde{a}[[x_ix^m][vv]]\tilde{b} + \sum_{1\leq l\leq z} \frac{\alpha_l}{2}\tilde{a}_l[[x_ix^m][vv]]\tilde{b}_l\\
&=&\frac{1}{2} \tilde{a} x_ix^m[vv]\tilde{b}- \frac{1}{2} \tilde{a}[vv]x_ix^m\tilde{b}+\sum_{\substack{s+t=m\\   0< s \leq m, 0\leq t < m}}\frac{\beta_{st}}{2} \tilde{a}(x^sx_ix^t [vv]-[vv]x^sx_ix^t)\tilde{b}\\
&&+ \sum_{1\leq l\leq z} \frac{\alpha_l}{2}\tilde{a}_l(x_ix^m[vv]-[vv]x_ix^m)\tilde{b}_l \\
&& + \sum_{\substack{1\leq l\leq z, s+t=m\\   0< s \leq m, 0\leq t < m}} \frac{\beta_{st}\alpha_l}{2}\tilde{a}_l (x^sx_ix^t [vv]-[vv]x^sx_ix^t)\tilde{b}_l\\
&=&\frac{1}{2} a[vv]b- \frac{1}{2} \tilde{a}[vv]x_ix^m\tilde{b}+\sum_{\substack{s+t=m\\   0< s \leq m, 0\leq t < m}}\frac{\beta_{st}}{2}( \tilde{a}x^sx_ix^t [vv]\tilde{b}- \tilde{a}[vv]x^sx_ix^t\tilde{b})\\
&&+\sum_{1\leq l\leq z} \frac{\alpha_l}{2}\tilde{a}_lx_ix^m[vv]\tilde{b}_l-\tilde{a}_l[vv]x_ix^m\tilde{b}_l \\
&& +  \sum_{\substack{1\leq l\leq z, s+t=m\\   0< s \leq m, 0\leq t < m}} \frac{\beta_{st}\alpha_l}{2}(\tilde{a}_lx^sx_ix^t [vv]\tilde{b}_l-\tilde{a}_l[vv]x^sx_ix^t\tilde{b}_l).
\end{eqnarray*}
It is easy to check that  $avvb$ is larger than the following words,
$$
\tilde{a}vvx_ix^m\tilde{b}, \tilde{a}x^sx_ix^tvv\tilde{b}, \tilde{a}vvx^sx_ix^t\tilde{b},
$$
$$
\tilde{a}_lx_ix^mvv\tilde{b}_l,\tilde{a}_lvvx_ix^m\tilde{b}_l,\tilde{a}_lx^sx_ix^t vv\tilde{b}_l,\tilde{a}_lvvx^sx_ix^t\tilde{b}_l.
$$

Case 2.  If $v\neq x$, then we set   $b=x^m\tilde{b}$, where $m\geq 0$ and the first letter of $\tilde{b}$ is not $x$.

Case 2-1.   If  $m>0$,  then
$$
[u]'_{vv}=[u]_{vvx^m}|_{[vvx^m]\mapsto [vv](\textbf{ad}x)^m}.
$$
By Theorem \ref{th2.15}, we  have
$$
[u]_{vvx^m}= a[vvx^m]\tilde{b} + \sum_{1\leq l\leq z} \alpha_l\tilde{a}_l[vvx^m]\tilde{b}_l.
$$
It follows that
$$
[u]'_{vv}= a([vv](\textbf{ad}x)^m)\tilde{b} + \sum_{1\leq l\leq z} \alpha_l\tilde{a}_l([vv](\textbf{ad}x)^m)\tilde{b}_l.
$$
By Lemma \ref{lem2.18}, we have
$$
[vv](\textbf{ad}x)^m=[vv]x^m+\sum_{\substack{i+j=m\\   0< i \leq m, 0\leq j < m}} \beta_{ij}x^i[vv]x^j,
$$
which gives that
\begin{eqnarray*}
[u]'_{vv}&=& a([vv](\textbf{ad}x)^m)\tilde{b} + \sum_{1\leq l\leq z} \alpha_l\tilde{a}_l([vv](\textbf{ad}x)^m)\tilde{b}_l\\
&=&  a[vv]x^m\tilde{b}+\sum_{\substack{i+j=m\\   0< i \leq m, 0\leq j < m}} \beta_{ij}ax^i[vv]x^j\tilde{b}\\
&&+ \sum_{1\leq l\leq z} \alpha_l\tilde{a}_l([vv]x^m+\sum_{\substack{i+j=m\\   0< i \leq m, 0\leq j < m}} \beta_{ij}x^i[vv]x^j)\tilde{b}_l\\
&=&  a[vv]b+\sum_{\substack{i+j=m\\   0< i \leq m, 0\leq j < m}} \beta_{ij}ax^i[vv]x^j\tilde{b}\\
&&+ \sum_{1\leq l\leq z} \alpha_l\tilde{a}_l[vv]x^m\tilde{b}_l+\sum_{\substack{1\leq l\leq  z,  i+j=m\\   0< i \leq m, 0\leq j < m}} \beta_{ij} \tilde{a}_lx^i[vv]x^j\tilde{b}_l
\end{eqnarray*}
Moreover,  we can see that
$$
avvb>ax^ivvx^j\tilde{b}, \tilde{a}_lvvx^m\tilde{b}_l, \tilde{a}_lx^ivvx^j\tilde{b}_l.
$$

Case 2-2.  If $m=0$, then   $[u]'_{vv}=[u']'_{v'v'}|_{x_ix^l\mapsto [x_ix^l]}$,  where  $u'=a'v'v'b'=x_{i_1}x^{l_1}x_{i_2}x^{l_2}\cdots x_{i_p}x^{l_p}$ and   $v'=x_{i_s}x^{l_s}x_{i_{s+1}}x^{l_{s+1}}\cdots x_{i_t}x^{l_t}$
are written down  in the new letters of the form $x_ix^l$. By the inductive hypothesis,
$$
[u']'_{v'v'}=\alpha a'[v'v']b'+\sum \alpha_h'a_h'[v'v']b_h',
$$
where $avvb=a'v'v'b'>a_h'v'v'b_h'=a_hvvb_h$.
It follows that
$$
[u]'_{vv}=(\alpha a'[v'v']b'+\sum \alpha_h'a_h'[v'v']b_h')|_{x_ix^l\mapsto [x_ix^l]}.
$$
Since  $[v'v']|_{x_ix^l\mapsto [x_ix^l]}=[vv]$ and  for each new letter $x_ix^l$,
$$
[x_ix^l]=x_ix^l+\sum_{\substack{s+t=l\\   0< s \leq l, 0\leq t < l}} \beta_{st} x^sx_ix^t, x_ix^l>x^sx_ix^t,
$$
we can obtain
$$
[u]_{vv}=\alpha a[vv]b + \sum \alpha_j a_j[vv]b_j,
$$
where for each $j$, $avvb> a_jvvb_j$. \hfill $\square$

\section{Gr\"{o}bner-Shirshov bases  theory for   operated   associative  superalgebras}

\subsection{Free operated   associative  superalgebras}

In this subsection, we will give a construction of  free operated   associative  superalgebra on a $\mathbb{Z}_2$-graded set $X=X_{\bar{0}}\cup X_{\bar{1}}$.
\begin{definition}  An operated  associative (resp. Lie) superalgebra is  pair $(A, P)$, where  $A=A_{\bar{0}}\oplus A_{\bar{1}}$ is an associative (resp. Lie) superalgebra and   $P: A\rightarrow A$ is an  even linear map on $A$.
\end{definition}

\begin{definition} Let $(A, P_A)$ and $(B, P_B)$ be two  operated associative (resp. Lie) superalgebras.   A homomorphism $f: A\rightarrow B$ from  associative (resp. Lie)    superalgebra $A$  to $B$ is said to be a homomorphism of   operated  associative (resp. Lie) superalgebra from  $(A, P_A)$ to $(B, P_B)$,   if  $P_B\circ f=f\circ P_A$.
\end{definition}

\begin{definition}Let $(A, P)$  be an operated  associative (resp. Lie)  superalgebra.    An   super-subalgebra  $B$ (resp.  super-ideal $I$) of  $A$ is called an operated super-subalgebra  (resp.   super-ideal) of  $(A, P)$ if $P(B)\subseteq B$  (resp.  $P(I)\subseteq I$).
\end{definition}

\begin{definition}  A free  operated  associative (resp. Lie)  superalgebra on a $\mathbb{Z}_2$-graded $X=X_{\bar{0}}\cup X_{\bar{1}}$  is  a pair $((F, P), \textbf{i})$  with $(F,P)$ an  operated associative (resp. Lie)  superalgebra and $\textbf{i} : X \rightarrow F$ an even operated map from $X$ to $F$ such that if $(L, P)$ is
any  operated (resp. Lie)    superalgebra and $\textbf{j} : X \rightarrow L$ is an  even operated  map, then there is a unique  homomorphism     of  operated associative (resp. Lie)   superalgebras  $\phi_\textbf{j}:  F\rightarrow L$ such that  $\phi_\textbf{j}\textbf{i}=\textbf{j}$.
\end{definition}

For any set $Y$,  we denote by $\textbf{P}(Y)= \{P(y)|y\in Y\}.$ Let $X$ be a set. We will define the sets
$X^{(P,n)}$ and  $\langle P;X;n\rangle$  by induction on $n$.

For $n=0$,  define $X^{(P,0)}=X$ and $\langle P;X;0\rangle=S(X^{(P,0)})$.

Assume that we have defined $X^{(P,n-1)}$ and  $\langle P;X;n-1\rangle=S(X^{(P,n-1)})$. Note that
$$
\textbf{P}(\langle P;X;n-1\rangle)= \{P(y)|y\in \langle P;X;n-1\rangle\}.
$$
Define
$$
X^{(P,n)}= X\cup \textbf{P}(\langle P;X;n-1\rangle), \ \langle P;X;n\rangle=S(X^{(P,n)}).
$$
It is easy to see  that   $X^{(P,n)}\subset X^{(P,n+1)}$ and $\langle P;X;n\rangle\subset \langle P;X;n+1\rangle$ for any $n\geq 0$.

We will denote by
$$
X^{(P)}= \lim_{n\rightarrow \infty} X^{(P,n)} =\bigcup_{n\geq 0}X^{(P,n)}
$$
$$
\langle P;X\rangle= \lim_{n\rightarrow \infty} \langle P;X;n\rangle =\bigcup_{n\geq 0}\langle P;X;n\rangle.
$$
The elements  of $\langle P;X\rangle$ are called operated  associative word on $X$.

Define   map $\tau_n: \langle P;X; n\rangle\rightarrow S(X)$  by induction on $n$.   If $n=0$, then  $\tau_0(w)=w$ for any $w\in \langle P;X; 0\rangle$. Assume that we have defined $\tau_{n-1}: \langle P;X; n-1\rangle\rightarrow S(X)$. For any $w=w_1w_2\cdots w_m\in \langle P;X; n\rangle$, where $w_i\in X^{(P, n)}$, we define
$$
\tau_n(w)=\tau_n(w_1)\tau_n(w_2)\cdots \tau_n(w_m),
$$
where
$\tau_n(w_i)=w_i$ if $w_i\in X$ and $\tau_n(w_i)=\tau_{n-1}(w_i')$ if $w_i=P(w_i')$.

Let  $X=X_{\bar{0}}\cup X_{\bar{1}}$ is  a $\mathbb{Z}_2$-graded set.  For any $w\in \langle P;X; n\rangle$, if $n=0$, then $w=x_1x_2\cdots x_n\in S(X)$, where each  $x_i\in X$.    Define the parity of $w$ by
$
|w|=\sum_{i=1}^n|x_i|,
$
where $|x|=\bar{i}$, if $x\in X_{\bar{i}}, \ \ i=0,1$.  If $n>0$, then we define the parity of $w$ by
$
|w|=|\tau_n(w)|.
$

By the parity of the words,  $X^{(P,n)}$ (resp. $\langle P;X;n\rangle$) has the  $\mathbb{Z}_2$-grading
$$
X^{(P,n)}=X^{(P,n)}_{\bar{0}}\cup X^{(P,n)}_{\bar{1}},\
(resp. \  \langle P;X;n\rangle=\langle P;X;n\rangle_{\bar{0}}\cup \langle P;X;n\rangle_{\bar{1}}),
$$
where for  $i=0,1$,
$$
X^{(P, n)}_{\bar{i}}=\{w| w\in X^{(P, n)}, |w|=\bar{i}\}
$$
$$
\langle P;X;n\rangle_{\bar{i}}=\{w| w\in \langle P;X;n\rangle, |w|=\bar{i}\}.
$$
For $i=0,1$, we can get
$$
X^{(P)}_{\bar{i}}=\{w|w\in X^{(P)}, |w|=\bar{i}\}=\bigcup_{n\geq 0}X^{(P,n)}_{\bar{i}}.
$$
$$
\langle P;X\rangle_{\bar{i}}=\{w|w\in \langle P;X\rangle, |w|=\bar{i}\}=\bigcup_{n\geq 0}{\langle P;X\rangle}_{\bar{i}}.
$$
It is worth pointing out that $\langle P;X\rangle=S(X^{(P)})$ and $\langle P;X\rangle_{\bar{i}}=S(X^{(P)})_{\bar{i}}, i=0,1.$

Note that  $\mathbf{AS}(X^{(P,n)})=k\langle P;X;n\rangle$ is a  free nonunitary associative superalgebra on the $\mathbb{Z}_2$-graded set  $X^{(P,n)}=X^{(P,n)}_{\bar{0}}\cup X^{(P,n)}_{\bar{1}}$. Moreover,
$$
\mathbf{AS}(X^{(P,n)})=\mathbf{AS}(X^{(P,n)})_{\bar{0}}\oplus \mathbf{AS}(X^{(P,n)})_{\bar{1}},
$$
where  $\mathbf{AS}(X^{(P,n)})_{\bar{0}}=k\langle P;X;n\rangle_{\bar{0}}$ and $\mathbf{AS}(X^{(P,n)})_{\bar{1}}=k\langle P;X;n\rangle_{\bar{1}}$.

Define
$$
\mathbf{OAS}(X)=k\langle P;X\rangle=\mathbf{OAS}(X)_{\bar{0}}\oplus \mathbf{OAS}(X)_{\bar{1}},
$$
where $\mathbf{OAS}(X)_{\bar{0}}=k\langle P;X\rangle_{\bar{0}}$ and   $\mathbf{OAS}(X)_{\bar{1}}=k\langle P;X\rangle_{\bar{1}}$.

It is easy to check that
$$
\mathbf{OAS}(X)=\bigcup_{n\geq 0}\mathbf{AS}(X^{(P,n)}),
$$
$$
\mathbf{OAS}(X)_{\bar{0}}=\bigcup_{n\geq 0}\mathbf{AS}(X^{(P,n)})_{\bar{1}},\ \mathbf{OAS}(X)_{\bar{1}}=\bigcup_{n\geq 0}\mathbf{AS}(X^{(P,n)})_{\bar{1}}.
$$

Define the even  linear map  $\mathscr{P}:\mathbf{OAS}(X)\rightarrow \mathbf{OAS}(X)$ by extending the even  map
 $$
 \mathscr{P}:\langle P;X\rangle\rightarrow \langle P;X\rangle, \mathscr{P}(w)=P(w), \forall\ w\in \textbf{OAW}(X).
 $$

\begin{theorem} (a)\   $(\mathbf{OAS}(X), \mathscr{P})$ is a nonunitary operated  associative superalgebra.

(b)\ $((\mathbf{OAS}(X), \mathscr{P}), \textbf{i}$)   is a   free nonunitary operated associative superalgebra on the $\mathbb{Z}_2$-graded set $X=X_{\bar{0}}\cup X_{\bar{1}}$, where $\textbf{i}: X\rightarrow \mathbf{OAS}(X)$ is the inclusion map.
\end{theorem}
{\bf Proof.}  The proof     is straightforward.  \hfill $\square$

 \subsection{Composition-Diamond lemma    for   operated   associative  superalgebras}

Let $\star$ be a symbol with   $\star\notin X$. By a
$\star$-word   we mean any expression in $\langle P; X\cup
\{\star\}\rangle$   with only one occurrence of $\star$. The set of all
  $\star$-word     on $X$ is denoted by $\langle P;  X \rangle^\star$.
Similarly, a  $(\star_1,\star_2)$-word is an element of  $\langle P; X\cup
\{\star_1, \star_2\}\rangle$ with only one occurrence of $\star_1$ and  one occurrence of $\star_2$. The set of all $(\star_1,\star_2)$-words is denoted by $\langle P;  X \rangle^{\star_1, \star_2}$.

Let $\pi \in \langle P;  X \rangle^\star$     and $s\in \mathbf{OAS}(X)$. Then we call
$
\pi|_{s}=\pi|_{\star\mapsto s}  \
$
an   $s$-word.

If   $u \in X^{(P)}=X \cup \textbf{P}( \langle P;X\rangle)$, where $\textbf{P}(\langle P;X\rangle)=\{P(w)|w\in \langle P;X\rangle\}$,  then $u$  is called   prime. For any $u\in \langle P;X\rangle$, $u$ can be expressed   uniquely in the canonical form
$
u=u_1u_2\cdots u_n,\ n\geq 1,
$
where each $u_i$ is prime. The number  $n$  is called the breath of $u$, which is    denoted  by $bre(u)$. The degree of $u$, denoted by $Deg(u)$,  is defined to be   the total number  of all occurrences of all $x\in X$ and $P$ in $u$.  For any     $u \in \langle P;X\rangle$, define   the depth of $u$ to be
$
dep(u)=\min\{n|u\in \langle P;X;n\rangle\}.
$

For any
$
u=u_1u_2\cdots u_n,\  n\geq 1,
$
where each $u_i$ is prime, let us denote by
$$
wt(u)=(Deg(u), bre(u), u_1,u_2,\cdots, u_n).
$$

Let $X$  be a  well-ordered set  with the order  $>_X$.    Define the  Deg-lex order  $>_{_{Dl}} $ on $\langle P;X\rangle$ as follows. For any $u=u_1u_2\cdots u_n, v=v_1v_2\cdots v_m$,   where    $u_i,  v_j$  are  prime, define
$$
u>_{_{Dl}}v \ \ \mbox{if}\ \ wt(u)>wt(v)\ \ \mbox{lexicographically},
$$
where if $u_i= P(\tilde{u}_{i} ), v_i=P(\tilde{v}_{i} )$ and $Deg(u_i)=Deg(v_i)$, then $u_i>_{_{Dl}}v_i$  if  $\tilde{u}_i>_{_{Dl}} \tilde{v}_i$  by induction.

For any $f\in \mathbf{OAS}(X)$, let $\overline{f}$ be the   leading  word  of $f$ with respect to the order $>_{_{Dl}}$ on $\langle P;  X\rangle$.
The coefficient of  $\overline{f}$ is  denoted by $lc(f)$.

Note that for any $S\subseteq  \mathbf{OAS}(X)$, there is  a homogenous  subset  $\widetilde{S}\subseteq \mathbf{OAS}(X)$ such that $\textbf{Id}(S)=\textbf{Id}(\widetilde{S})$, where $\textbf{Id}(S)$ is the operated superideal of $\mathbf{OAS}(X)$ generated by $S$.

\begin{definition}
Let $f, g \in \mathbf{OAS}(X)$ be  homogenous.  Define  two kinds of compositions as follows.
\begin{enumerate}
\item[(a)]\
If  $w=\bar{f}a=b\bar{g}$ with  $bre(w)<bre(\bar{f})+bre(\bar{g})$, then
$$
\langle f,g\rangle_w= lc(f)^{-1}fa- lc(g)^{-1}bg
$$
is called the intersection composition of $f$ and $g$ with respect to the ambiguity $w$.

\item[(b)]\  If $w=\overline{f}=\pi|_{\overline{g}}$, where $\pi\in \langle P;  X\rangle^\star$,    then
$$
\langle f,g\rangle_w=lc(f)^{-1}f-lc(g)^{-1}\pi|_{g}
$$
is called the inclusion composition of $f$ and $g$ with respect to  the ambiguity  $w$.
\end{enumerate}
\end{definition}

If $S$ is a homogenous  subset of $\mathbf{OAS}(X)$ and  $h\in \mathbf{OAS}(X)$, then $h$ is called trivial modulo $(S, \leq_{_{Dl}})$ if
$$
h=\sum \alpha_i  \pi_i|_{s_i}
$$
where each $\alpha_i\in k,
  \ s_i\in S$,   $\pi_i|_{_{\overline{s_i}}}\leq_{_{Dl}} \overline{h}$. If
this is the case, then we write
$$
h\equiv 0\  mod(S, \leq_{_{Dl}}).
$$

\begin{definition} Let $S$ be a homogenous subset of $\mathbf{OAS}(X)$.
A homogenous set $S$ is called a
Gr\"{o}bner-Shirshov basis  in   $\mathbf{OAS}(X)$ if any
composition $\langle f,g\rangle_w$ of  $f, g\in S$ is trivial modulo $(S, \leq_{_{Dl}})$.
\end{definition}

\begin{theorem}\label{cdla}{\em(Composition-Diamond lemma for operated associative superalgebra)}\ \  Let $S$ be a homogenous subset of $\mathbf{OAS}(X)$
 and $\textbf{Id}(S)$ the operated superideal of
 $\mathbf{OAS}(X)$ generated by $S$.  Then the following
statements are equivalent:
 \begin{enumerate}
\item[(a)] The set  $S $ is a Gr\"{o}bner-Shirshov basis in $\mathbf{OAS}(X)$.
\item[(b)] For any $ f\in \textbf{Id}(S)$,  we have $\bar{f}=\pi|_{\overline{s}}$
for some $\pi \in \langle P;  X\rangle^\star$ and $s\in S$.
 \item[(c)] The set
 $$
 Irr(S) = \{ w\in  \langle P;  X\rangle|  w \neq
\pi|_{\overline{s}}
 \mbox{ for  any} \ \pi\in \langle P;  X\rangle^\star \ \mbox{and } s\in S\}
 $$
is a $k$-basis of $\mathbf{OAS}(X|S)=\mathbf{OAS}(X)/\textbf{Id}(S)$.
\end{enumerate}
\end{theorem}
{\bf Proof. }  The proof is similar to the   Composition-Diamond lemma for operated  associative algebra in \cite{bcq2010}.  \hfill $\square$

\section{Free    operated   Lie   superalgebras}

\subsection{Operated Super-Lyndon-Shirshov  words and monomials}

In this subsection, we   define the  operated super-Lyndon-Shirshov (SLS) words and monomials on a  $\mathbb{Z}_2$-graded set $X=X_{\bar{0}}\cup X_{\bar{1}}$, which are   generalizations of the classical    super-Lyndon-Shirshov words and  monomials.

Let $X=X_{\bar{0}}\cup X_{\bar{1}}$ be a  $\mathbb{Z}_2$-graded set. We will define the set $\{P;X;n\}$  by induction on $n$.

For $n=0$,  define   $\{P;X;0\}= X^{**}$, the set of all nonassociative words on $X$.

Assume that we have defined   $\{P;X;n-1\}$. Let us denote by
$$
\textbf{P}(\{P;X;n-1\})= \{P(y)|y\in \{P;X;n-1\}\}.
$$
Let
$$
\{ P;X;n\}=(X\cup \textbf{P}(\{P;X;n-1\}))^{**}.
$$
It is easy to see  that    $\{ P;X;n\}\subset \{P;X;n+1\}$ for any $n\geq 0$.  Let
$$
\{ P;X\}=\bigcup_{n\geq 0} \{ P;X;n\}.
$$
The elements  of $\{ P;X\}$ are called operated nonassociative words on $X.$

For any nonassociative word $(u)\in \{ P;X\}$, define the parity of $u$ by $|(u)|=|u|$. For example, if $(u)=(P(x_0)(x_1P(x_1)))$, where $x_0\in X_{\bar{0}}$ and $x_{\bar{1}}\in X_{\bar{1}}$, then
$$
|(u)|=|(P(x_0)(x_1P(x_1)))|=|P(x_0)x_1P(x_1)|=\bar{0}.
$$
Thus, $\{ P;X\}$ is a $\mathbb{Z}_2$-graded set, i.e.
$$
\{ P;X\}=\{P;X\}_{\bar{0}}\cup  \{ P;X\}_{\bar{1}},
$$
where $\{P;X\}_{\bar{i}}=\{(w)|(w)\in \{P;X\}, |(w)|=\bar{i}\}, i=0,1$. Moreover, for each $n\geq 0$

$$
\{ P;X;n\}=\{P;X;n\}_{\bar{0}}\cup  \{ P;X;n\}_{\bar{1}},
$$
where $\{P;X;n\}_{\bar{i}}=\{(w)|(w)\in \{P;X;n\}, |(w)|=\bar{i}\}, i=0,1$. We also define the depth of $(w)$ by $dep((w))=dep(w)$ for any $(w)\in \{ P;X\}$.

We will define the $\mathbb{Z}_2$-graded sets
$$
\mathfrak{X}^{(P,n)}=\mathfrak{X}^{(P,n)}_{\bar{0}}\cup  \mathfrak{X}^{(P,n)}_{\bar{1}}
$$
$$
\textbf{OSLSW}(X, n)= \textbf{OSLSW}(X, n)_{\bar{0}}\cup \textbf{OSLSW}(X, n)_{\bar{1}}
$$
by induction on $n$.

For $n=0$, set $\mathfrak{X}^{(P,0)}=\mathfrak{X}^{(P,0)}_{\bar{0}}\cup  \mathfrak{X}^{(P,0)}_{\bar{1}}$, where $\mathfrak{X}^{(P,0)}_{\bar{0}}=X_{\bar{0}},    \mathfrak{X}^{(P,0)}_{\bar{1}}=X_{\bar{1}}$.  Let  $SLSW(\mathfrak{X}^{(P,0)})$  be the set of all   super-Lyndon-Shirshov   words on   the $\mathbb{Z}_2$-graded set $\mathfrak{X}^{(P,0)}=\mathfrak{X}^{(P,0)}_{\bar{0}}\cup\mathfrak{X}^{(P,0)}_{\bar{1}}$. Define
$$
\textbf{OSLSW}(X, 0)= \textbf{OSLSW}(X, 0)_{\bar{0}}\cup \textbf{OSLSW}(X, 0)_{\bar{1}}
$$
where  for $i=0,1$,
$
\textbf{OSLSW}(X, 0)_{\bar{i}}=\{w\in SLSW(\mathfrak{X}^{(P,0)})||w|=\bar{i}\}.
$

Assume that we have defined the $\mathbb{Z}_2$-graded set
$$
\mathfrak{X}^{(P,n-1)}=\mathfrak{X}^{(P,n-1)}_{\bar{0}}\cup  \mathfrak{X}^{(P,n-1)}_{\bar{1}}
$$
and
$$
\textbf{OSLSW}(X, n-1)= \textbf{OSLSW}(X, n-1)_{\bar{0}}\cup \textbf{OSLSW}(X, n-1)_{\bar{1}},
$$
where  for $i=0,1$,
$
\textbf{OSLSW}(X, n-1)_{\bar{i}}=\{w\in SLSW(\mathfrak{X}^{(P,n-1)})||w|=\bar{i}\}.
$

Let us denote by
$$
\textbf{P}(\textbf{OSLSW}(X, n-1))=\textbf{P}(\textbf{OSLSW}(X, n-1))_{\bar{0}} \cup \textbf{P}(\textbf{OSLSW}(X, n-1)))_{\bar{1}}.
$$
where $\textbf{P}(\textbf{OSLSW}(X, n-1))_{\bar{i}}=\textbf{P}(\textbf{OSLSW}(X, n-1))_{\bar{i}}), i=0,1$.
Define
$$
\mathfrak{X}^{(P,n)}=\mathfrak{X}^{(P,n)}_{\bar{0}}\cup  \mathfrak{X}^{(P,n)}_{\bar{1}},
$$
where  $\mathfrak{X}^{(P,n)}_{\bar{i}}=X_{\bar{i}} \cup \textbf{P}(\textbf{OSLSW}(X, n-1))_{\bar{i}}, i=0,1$.

Let $\succ$ be the restriction of $>_{_{Dl}}$  on $\mathfrak{X}^{(P,n)}$ and $SLSW(\mathfrak{X}^{(P,n)})$  be the set of all super-Lyndon-Shirshov   words  on the $\mathbb{Z}_2$-graded set
$
\mathfrak{X}^{(P,n-1)}=\mathfrak{X}^{(P,n-1)}_{\bar{0}}\cup  \mathfrak{X}^{(P,n-1)}_{\bar{1}}
$
with respect to the lex-order $\succ_{lex}$.    Define
$$
\textbf{OSLSW}(X, n)= \textbf{OSLSW}(X, n)_{\bar{0}}\cup \textbf{OSLSW}(X, n)_{\bar{1}},
$$
where  for $i=0,1$,  $\textbf{OSLSW}(X, n)_{\bar{i}}=\{w\in SLSW(\mathfrak{X}^{(P,n)})||w|=\bar{i}\}.$

Define
$$
\textbf{OSLSW}(X)=\bigcup_{n\geq 0} \textbf{OSLSW}(X, n),
$$
$$
\textbf{OSLSW}(X)_{\bar{i}}=\bigcup_{n\geq 0} \textbf{OSLSW}(X, n)_{\bar{i}}, i=0,1.
$$
It is easily seen that  $\textbf{OSLSW}(X)$ has a  $\mathbb{Z}_2$-grading
$$
\textbf{OSLSW}(X)=\textbf{OSLSW}(X)_{\bar{0}} \cup  \textbf{OSLSW}(X)_{\bar{1}}.
$$
The elements of  $\textbf{OSLSW}(X)$ are called operated super-Lyndon-Shirshov   words on the $\mathbb{Z}_2$-graded set
$
X=X_{\bar{0}}\cup  X_{\bar{1}}.
$
Moreover,  an operated super-Lyndon-Shirshov   word $u$ is said to be pure, if $u\neq vv$ with $|v|=1$.

Define the bracketing map $[\ ]:\textbf{OSLSW}(X, n)\rightarrow  \{ P;X\}$ by induction on $n$ as follows.

For any  $w\in \textbf{OSLSW}(X, 0)$,  define $[w]$   the  unique super-Lyndon-Shirshov   monomial  related to the super-Lyndon-Shirshov word $w$.

Assume that we have defined  the bracketing map
$$
[\ ]:\textbf{OSLSW}(X, n-1)\rightarrow  \{ P;X\}.
$$

For  $w\in \mathfrak{X}^{(P,n)}$, if  $w\in X$, then we   define $[w]=w$. Otherwise,  $w\in \textbf{P}(\textbf{OSLSW}(X, n-1))$, i.e.  $w=P(\tilde{w})$, where $\tilde{w}\in \textbf{OSLSW}(X, n-1)$, and we  define $[w]=P([\tilde{w}])$, where $[\tilde{w}]$ is defined by induction.  Set
$$
[\mathfrak{X}^{(P,n)}]=\{[w]| w\in \mathfrak{X}^{(P,n)}\}.
$$
Then we have
$$
[\mathfrak{X}^{(P,n)}]=[\mathfrak{X}^{(P,n)}]_{\bar{0}}\cup [\mathfrak{X}^{(P,n)}]_{\bar{1}},
$$
where
$$
[\mathfrak{X}^{(P,n)}]_{\bar{0}}=\{[w]| w\in \mathfrak{X}^{(P,n)}_{\bar{0}}\}, [\mathfrak{X}^{(P,n)}]_{\bar{1}}=\{[w]| w\in \mathfrak{X}^{(P,n)}_{\bar{1}}\}.
$$
The order $\succ$ on $\mathfrak{X}^{(P,n)}$  induces an order (still denoted  by $\succ$) on  $[\mathfrak{X}^{(P,n)}]$   by $[u] \succ[ v ]$  if $u\succ v$ for any $u, v\in  \mathfrak{X}^{(P,n)}$.

If $u=u_1u_2\cdots u_m\in \textbf{OSLSW}(X, n)$, where each  $u_i\in \mathfrak{X}^{(P,n)}$, then we define $[u]$ by two steps. Firstly, we define  $[u]^{(1)}=[ u_1 ][ u_2]\cdots [ u_m ]$.  Note that $[ u_1 ][ u_2]\cdots [ u_m ]$ is a  super-Lyndon-Shirshov  word on $[\mathfrak{X}^{(P,n)}]=[\mathfrak{X}^{(P,n)}_{\bar{0}}]  \cup [\mathfrak{X}^{(P,n)}_{\bar{1}}]$.  Secondly, we define $[[u]^{(1)}]^{(2)}=[[ u_1 ][ u_2]\cdots [ u_m ]]$, the  unique super-Lyndon-Shirshov   monomial  related to the super-Lyndon-Shirshov word $[ u_1 ][ u_2]\cdots [ u_m ]$. Set
$$
[u]  =[[u]^{(1)}]^{(2)}=[[ u_1 ][ u_2]\cdots [ u_m ]].
$$
We will denote by
$$
\textbf{OSLSM}(X)=\{[w]|w\in \textbf{OSLSW}(X)\}.
$$
It is   obvious that
$$
\textbf{OSLSM}(X)=\textbf{OSLSM}(X)_{\bar{0}}\cup \textbf{OSLSM}(X)_{\bar{1}},
$$
where $\textbf{OSLSM}(X)_{\bar{i}}=\{[w]\in \textbf{OSLSW}(X)| |[w]|=\bar{i}\}, i=0,1.$ The elements of $\textbf{OSLSM}(X)$  are called  operated super-Lyndon-Shirshov  monomials.

\subsection{Free   Lie  operated superalgebras}

In this subsection, we prove that the  set    of  all operated super-Lyndon-Shirshov   monomials  $\textbf{OSLSM}(X)$ is  a linear basis of the  free  operated Lie superalgebra on the $\mathbb{Z}_2$-graded set     $X=X_{\bar{0}}\cup X_{\bar{1}}$.

Note that  $(\mathbf{OAS}(X), \mathscr{P})$ is a   free operated nonunitary associative superalgebra on the $\mathbb{Z}_2$-graded set $X=X_{\bar{0}}\cup X_{\bar{1}}$, where $\mathbf{OAS}(X)=\mathbf{OAS}(X)_{\bar{0}}\oplus \mathbf{OAS}(X)_{\bar{1}}$.  It follows that
$(\mathbf{OAS}(X), [, ] \mathscr{P})$ is an  operated   Lie superalgebra under the  superbracket
$$
[uv]=uv-(-1)^{|u||v|}vu,
$$
where $u$ and $v$ are homogenous elements.

Let  $\textbf{OLS}(X)$  be the operated Lie  super-subalgebra of   $(\mathbf{OAS}(X), [, ],  \mathscr{P})$     generated by the $\mathbb{Z}_2$-graded set  $X=X_{\bar{0}}\cup X_{\bar{1}}$.

\begin{lemma}\label{lem2.6}
For any  $(u)\in  \{P;X\}$,   $(u)$ has a representation
$$
(u)=\sum \alpha_i[u_i],
$$
in $\textbf{OLS}(X)$, where each $\alpha_i\in k$ and $[u_i]\in \textbf{OSLSM}(X)$.
\end{lemma}
{\bf Proof.} The proof is by induction on the depth of  $(u)$. If $dep((u))=0$, then $(u)\in X^{**}$. By Lemma \ref{le2.2},
$
(u)=\sum \alpha_i[u_i],
$
where each  $a_i\in k, [u_i]\in \textbf{SLSM}(X)\subseteq \textbf{OSLSM}(X)$. Suppose that the result is true for any $(u)\in \{P;X\}$ with $dep((u))<n$. Let $(u)\in \{P;X\}$ with $dep((u))=n$. Assume that $(u)=(w_1w_2\cdots w_m)$, where each $w_i\in X\cup \textbf{P}(\{P;X;n-1\})$.

Case 1.  If $m=1$, then $(u)=P((\tilde{u}))$ with  $dep((\tilde{u}))=n-1$.  By the induction hypothesis,  we have
$
(\tilde{u})=\sum \alpha_i[\tilde{u}_i],
$
where each $[\tilde{u}_i]\in \textbf{OSLSM}(X)$.  It follows that
$
(\tilde{u})=P((\tilde{u}))=\sum \alpha_iP([\tilde{u}_i]),
$
where each $P([\tilde{u}_i])\in \textbf{OSLSM}(X)$.

Case 2.  If $m\geq 1$, then $dep((w_i))\leq n$.  By induction   hypothesis and Case 1, we can assume that $(\tilde{u})=([w_1][w_2]\cdots [w_m])$, where each $[w_j]\in [\mathfrak{X}^{(P,n)}]$. Note that  $(\tilde{u})=([w_1][w_2]\cdots [w_m])$ is a   nonassociative word on the   letters $[w_1],[w_2],\cdots, [w_m]\in [\mathfrak{X}^{(P,n)}]=[\mathfrak{X}^{(P,n)}]_{\bar{0}}\cup [\mathfrak{X}^{(P,n)}]_{\bar{1}}$. By Lemma \ref{le2.2}, we have
$$
(u)=([w_1][w_2]\cdots [w_m])=\sum \alpha_i[u_i],
$$
where each $[u_i]\in SLSM([w_1],[w_2],\cdots, [w_m])\subseteq \textbf{OSLSM}(X).$ Therefore, the result is true for any $(u)\in  \{P;X\}$.   \hfill $\square$

\begin{lemma}\label{lem2.7}
If  $u\in \textbf{OSLSW}(X)$, then $\overline{[u]}=u$ in $\textbf{OLS}(X)\subseteq \mathbf{OAS}(X)$ with respect to  the order $>_{_{Dl}}$.
\end{lemma}
{\bf Proof.}  The proof is by induction on the depth of  $u$. If $dep(u)=0$, then $u\in \textbf{OSLSW}(X, 0)$. By Lemma \ref{le2.2}, we can get $\overline{[u]}=u$. Suppose the result is true for any $u\in \textbf{OSLSW}(X)$ with $dep(u)<n$.  Let $dep(u)=n$ and $u=w_1w_2\cdots w_m$, where each $w_i\in X\cup \textbf{P}(\textbf{OSLSW}(X)).$

Case 1.  If $m=1$,  then $u=P(\tilde{u})$, where $\tilde{u}\in  \textbf{OSLSW}(X)$  and  $dep(\tilde{u})=n-1$.  By the induction hypothesis,  we have $\overline{[\tilde{u}]}=\tilde{u}$, which gives that $\overline{[u]}=P(\overline{[\tilde{u}]})=P(\tilde{u})=u$.

Case 2.  If $m>1$,  then $dep(w_i)\leq n$ for $1\leq i\leq m$.  By the induction hypothesis and Case 1, we have $\overline{[w_i]}=w_i$. By the definition of $[u]$, we have
$[u]=[[w_1][w_2]\cdots [w_m]]$. It follows that
$$
\overline{[u]}= \overline{[w_1][w_2]\cdots [w_m]} =w_1w_2\cdots w_m=u.
$$
Therefore, we can obtain the  result  by induction.  \hfill $\square$

\begin{theorem}\label{th4.3} The set $\textbf{OSLSM}(X)$     is  a linear basis of the   operated Lie  superalgebra $\textbf{OLS}(X)$.
\end{theorem}
{\bf Proof.}  Firstly, we shows that $\textbf{OSLSM}(X)$ is  a linear independent set. That is to show that for any distinct elements
$[u_1], [u_2],\cdots, [u_m]\in \textbf{OSLSM}(X)$, $[u_1], [u_2],\cdots, [u_m] $ are  linear independent.  We can prove it by induction on $m$. The results hold trivially when $m=1$. Assume the result is true  for any $m-1$ distinct elements in $\textbf{OSLSM}(X)$. Without loss generality,  we can assume that     $[u_1]>_{Dl}[u_2]>_{Dl}\cdots>_{Dl}[u_m]$.  Consider the scalars $k_1, k_2,\cdots, k_m$ such that
$$
k_1[u_1]+k_2[u_2]+\cdots+k_m[u_m]=0.
$$
By Lemma \ref{lem2.7}, $\overline{[u_i]}=u_i$ for $i=1,2,\cdots, m$.   If $k_1\neq 0$, then
$$
\overline{k_1[u_1]+k_2[u_2]+\cdots+k_m[u_m]}=u_1,
$$
a contradiction. This gives that $k_1=0$. Thus we have
$$
k_2[u_2]+k_3[u_3]+\cdots+k_m[u_m]=0.
$$
Applying the induction hypothesis, $[u_2], [u_3],\cdots, [u_m]$ are  linear independent. It follows that  $k_2=k_3=\cdots=k_m=0$. Therefore,   $[u_1], [u_2],\cdots, [u_m] $ are  linear independent.

  Secondly, by Lemma \ref{lem2.6}, we can obtain any element  in  $\textbf{OLS}(X)$ can be linearly  expressed by  the elements  in  $\textbf{OSLSM}(X)$.  Therefore,   $\textbf{OSLSM}(X)$ is  a linear basis of the   operated Lie  superalgebra $\textbf{OLS}(X)$.\hfill $\square$

\begin{theorem}
$(\textbf{OLS}(X),\mathscr{P}) $ is a free  operated  Lie  superalgebra on the   $\mathbb{Z}_2$-graded set $X=X_{\bar{0}}\cup X_{\bar{1}}$ with a   linear basis $\textbf{OSLSM}(X)$.
\end{theorem}
{\bf Proof.}   By  Theorem \ref{th4.3},  we can obtain the result. \hfill $\square$

\section{Gr\"{o}bner-Shirshov bases  theory  for  operated  Lie  superalgebras}

\subsection{Properties of operated  super-Lyndon-Shirshov  words and monomials}

In this subsection, we   give some properties of  operated   super-Lyndon-Shirshov  words and monomials.

\begin{proposition}\label{pro5.1}  Let $w_1, w_2$  be  operated super-Lyndon-Shirshov  words and $w_1=e_1e_2, w_2=e_2e_3$  where $e_1, e_2, e_3$ are nonempty word. Then we have the following three results.

(1)   If $w_1=u$,  $w_2=v^n (n=1, 2)$, $u, v$ are  pure operated   super-Lyndon-Shirshov  word   and $u\neq v$,  then $e_1e_2e_3$ is a pure operated super-Lyndon-Shirshov  word.

(2)   If $w_1=uu$,   $w_2=v^n(n=1, 2)$, $u, v$ are  pure  operated    super-Lyndon-Shirshov  word   and $u\neq v$, then $e_1=u\tilde{e}_1$ and $e_1e_2e_3$ is a pure  operated super-Lyndon-Shirshov   word, where $\tilde{e}_1$ may be empty.

(3)   If $w_1=uu$,    $w_2=uu$, $u$  is a  pure  operated  super-Lyndon-Shirshov  word with $|u|=1$,  then $e_1=e_2=e_3=u$.
\end{proposition}
{\bf Proof.}  Let $n=\max\{dep(w_1), dep(w_2)\}$. Then $w_1, w_2$  are super-Lyndon-Shirshov words on the the $\mathbb{Z}_2$-graded set $\mathfrak{X}^{(P,n)}=\mathfrak{X}^{(P,n)}_{\bar{0}}\cup  \mathfrak{X}^{(P,n)}_{\bar{1}}.$   By  Proposition  \ref{pro2.10},  we can obtain the result.    \hfill $\square$

\begin{proposition}  \label{lem5.2}  Let $w_1, w_2=vv$  be   operated super-Lyndon-Shirshov words with $|v|=1$.  If  $w_1$ is a proper  subword of $w_2$, then  $w_2=\pi|_{w_1}v$ or $w_2=v\pi|_{w_1}$, where $v=\pi|_{w_1}$.
\end{proposition}
{\bf Proof.}   Let $n=dep(w_2)$. Then  $w_1, w_2=vv$ are super-Lyndon-Shirshov words on the the $\mathbb{Z}_2$-graded set $\mathfrak{X}^{(P,n)}=\mathfrak{X}^{(P,n)}_{\bar{0}}\cup  \mathfrak{X}^{(P,n)}_{\bar{1}}.$   By Proposition \ref{pro2.11}, we can get  $w_1$ is a subword of $v$. Let $v=v_1v_2\cdots v_m,$ where each $v_i\in \mathfrak{X}^{(P,n)}$. If $w_1$ is  a subword of some $v_i$, then $w_1=v_i$  or $v_i=\pi_i|_{w_i}$.  Otherwise, we have  $v=aw_1b$. In either case,  we can get the result.
 \hfill $\square$

\begin{theorem}  \label{th5.3} {\em ( \cite{qc2017})} \label{th5.3} Let  $\pi\in \langle P;  X\rangle^\star$ and $v, \pi|_v\in \textbf{OSLSM}(X)$ are  pure operated super-Lyndon-Shirshov words.  Then there  is a  $\pi_{_1}\in \langle P;  X\rangle^\star$ and $c\in \langle P;  X\rangle$ such that
$
[ \pi|_v] =[\pi_{_1}|_{[vc]}],
$
where $c$ may be empty. Define
$$
\llbracket\pi|_v\rrbracket_{v}=[\pi_{_1}|_{[vc]}]|_{[vc]\mapsto [[\cdots[[[v][c_1]][c_2]]\cdots ][c_m]]},
$$
where   $c=c_1c_2\cdots c_m$ with  each  $c_i\in \textbf{OSLSM}(X)$,  $c_{1} \preceq_{{\rm lex}} c_{2}\preceq_{{\rm lex}} \cdots \preceq_{{\rm lex}} c_m$.   Then we have
$$
\llbracket\pi|_v\rrbracket_{v}=\pi|_{[v]}+\sum \alpha_i \pi_i|_{[v]},
$$
where each $\alpha_i\in k$ and $\pi_i|_{v}<_{_{\rm Dl}}\pi|_{v}$. It follows that $ \overline{\llbracket\pi|_v\rrbracket_{v}}=\pi|_v$   with respect to the order $>_{_{\rm Dl}}$.
\end{theorem}

  Let  $\pi\in \langle P;  X\rangle^\star$ and $v, \pi|_{vv}\in \textbf{OSLSW}(X)$ with $|v|=1$.  We will define $\llbracket\pi|_{vv}\rrbracket'_{vv}$ on $ \pi|_{vv}$
  by induction on the  depth of $\pi|_{vv}$ as follows.

  If $dep(\pi|_{vv})=0$, then  $\pi|_{vv}\in  \textbf{OSLSW}(X, 0)$ and we define $\llbracket\pi|_{vv}\rrbracket'_{vv}=[\pi|_{vv}]'_{vv}$.

   Assume that we have defined $\llbracket\pi|_{vv}\rrbracket'_{vv}$ for $dep(\pi|_{vv})<n$.    Let $\pi|_{vv}=w_1w_2\cdots w_m$, where each $w_i\in X\cup \textbf{P}(\textbf{OSLSW}(X))$ and $m\geq 1$.

   Case 1. If $m=1$, then $\pi|_{vv}=P(\pi_1|_{vv})$ with  $dep(\pi_1|_{vv})=n-1$.  Define
   $$
   \llbracket\pi|_{vv}\rrbracket'_{vv}=P(\llbracket\pi_1|_{vv}\rrbracket'_{vv}),
   $$
   where $\llbracket\pi_1|_{vv}\rrbracket'_{vv}$  is defined by induction.

   Case 2.  If $m>1$, there are two case to consider.

   Case 2-1.
   If $vv$ is located  in  the word   $w_i$, i.e.  $w_i=\pi_i|_{vv}$, then we define
   $$
   \llbracket\pi|_{vv}\rrbracket'_{vv}= \llbracket\tilde{\pi}|_{w_i}\rrbracket_{w_i}|_{[w_i]\mapsto \llbracket\pi_i|_{vv}\rrbracket'_{vv}},
   $$
   where $\tilde{\pi}|_{w_i}=\pi|_{vv}$ and $\llbracket\pi_i|_{vv}\rrbracket'_{vv}$ is defined by the  Case 1.

    Case 2-2.  If $\pi|_{vv}=w_1\cdots w_{t-1}vvw_{s+1}\cdots w_{m}$ and  $v=w_tw_{t+1}\cdots w_{s}, s>1$, where each $w_i\in X\cup \textbf{P}(\textbf{OSLSW}(X))$.    Note that  $[v]^{(1)}=[w_t][w_{t+1}]\cdots [w_{s}]$  and $[\pi|_{vv}]^{(1)} =[w_1][w_2]\cdots [w_m]$ are  super-Lyndon-Shirshov word on letters $[w_1],[w_2,]\cdots, [w_m]$ with $|[v]^{(1)}|=1$.  Define
   $$
    \llbracket\pi|_{vv}\rrbracket'_{vv}=[[\pi|_{vv}]^{(1)}]'_{[v]^{(1)}[v]^{(1)}}.
   $$

\begin{theorem}  \label{th5.4}  Let  $\pi\in \langle P;  X\rangle^\star$ and $v, \pi|_{vv}\in \textbf{OSLSW}(X)$ be  pure operated super-Lyndon-Shirshov words with $|v|=1$.      Then we have
$$
\llbracket\pi|_{vv}\rrbracket_{vv}'=\alpha \pi|_{[vv]}+\sum \alpha_i \pi_i|_{[vv]},
$$
where each $\alpha_i\in k$ and $\pi_i|_{vv}<_{_{\rm Dl}}\pi|_{vv}$. It follows that $ \overline{\llbracket\pi|_{vv}\rrbracket'_{vv}}=\pi|_{vv}$   with respect to the order $>_{_{\rm Dl}}$.
\end{theorem}
{\bf Proof.} The proof is by induction on the depth of  $\pi|_{vv}$.  If $dep(\pi|_{vv})=0$, then   $\pi|_{vv}\in  \textbf{OSLSW}(X, 0)$. It follows that    the result is true by  Theorem \ref{th2.15}.

Assume  that the result is true  for  the word $\pi|_{vv}$ with $dep(\pi|_{vv})<n$. We will prove it for $\pi|_{vv}$ with $dep(\pi|_{vv})=n$.

Let  $\pi|_{vv}=w_1w_2\cdots w_m$, where each $w_i\in X\cup \textbf{P}(\textbf{OSLSW}(X))$ and $m\geq 1$.

Case 1. If $m=1$, then   $\pi|_{vv}=P(\pi_1|_{vv})$ with  $dep(\pi_1|_{vv})=n-1$ and
$$
\llbracket\pi|_{vv}\rrbracket'_{vv}=P(\llbracket\pi_1|_{vv}\rrbracket'_{vv}).
$$
Since $dep(\pi_1|_{vv})=n-1$, by the induction hypothesis, we have
$$
\llbracket\pi_1|_{vv}\rrbracket_{vv}'=\alpha \pi_1|_{[vv]}+\sum \alpha_i \pi_{1i}|_{[vv]},
$$
where each $\pi_{1i}|_{vv} \leq \pi_1|_{vv}$.  It follows that
\begin{eqnarray*}
\llbracket\pi|_{vv}\rrbracket'_{vv}&=&P(\llbracket\pi_1|_{vv}\rrbracket'_{vv})\\
&=&P(\alpha \pi_1|_{[vv]}+\sum \alpha_i \pi_{1i}|_{[vv]})\\
&=&\alpha P(\pi_1|_{[vv]})+\sum \alpha_i P(\pi_{1i}|_{[vv]})\\
&=&\alpha \pi|_{[vv]}+\sum \alpha_i \pi_{i}|_{[vv]},
\end{eqnarray*}
where  $\pi_{i}=P(\pi_{1i})$.  Since  $\pi_1|_{vv}> \pi_{1i}|_{vv}$, we have  $\pi|_{vv}=P(\pi_{1}|_{vv})>  P(\pi_{1i}|_{vv})=\pi_{i}|_{vv}$.

Case 2.  If $m>1$ and $vv$ is located  in  the word   $w_i$, then
$$
   \llbracket\pi|_{vv}\rrbracket'_{vv}= \llbracket\tilde{\pi}|_{w_i}\rrbracket_{w_i}|_{[w_i]\mapsto \llbracket\pi_i|_{vv}\rrbracket'_{vv}},
$$
Since $w_i$ is a pure operated super-Lyndon-Shirshov word,  by Theorem \ref{th5.3},
$$
\llbracket\tilde{\pi}|_{w_i}\rrbracket_{w_i}=\tilde{\pi}|_{[w_i]}+\sum \alpha_i \tilde{\pi}_{i}|_{[w_i]},
$$
where $\tilde{\pi}|_{w_i}> \tilde{\pi}_{i}|_{w_i}$.  Moreover, by the Case 1,
$$
\llbracket\pi_i|_{vv}\rrbracket'_{vv}= \pi_i|_{[vv]}+\sum \alpha_{it} \pi_{it}|_{[vv]},
$$
where $\pi_i|_{vv}> \pi_{it}|_{vv}$.  Thus, we have
\begin{eqnarray*}
&& \llbracket\pi|_{vv}\rrbracket'_{vv}\\
&=&\llbracket\tilde{\pi}|_{w_i}\rrbracket_{w_i}|_{[w_i]\mapsto \llbracket\pi_i|_{vv}\rrbracket'_{vv}}\\
&=& ( \tilde{\pi}|_{[w_i]}+\sum \alpha_i \tilde{\pi}_{i}|_{[w_i]})|_{[w_i]\mapsto \llbracket\pi_i|_{vv}\rrbracket'_{vv}}\\
&=& \tilde{\pi}|_{\pi_i|_{[vv]}+\sum \alpha_{it} \pi_{it}|_{[vv]}}+\sum \alpha_i \tilde{\pi}_{i}|_{\pi_i|_{[vv]}+\sum \alpha_{it} \pi_{it}|_{[vv]}}\\
&=&   \pi|_{[vv]} +\sum\alpha_{it} \tilde{\pi}|_{\pi_{it}|_{[vv]}}+\sum \alpha_i \tilde{\pi}_{i}|_{\pi_i|_{[vv]}}+\sum  \alpha_i \alpha_{it} \tilde{\pi}_{i}|_{  \pi_{it}|_{[vv]}}.
\end{eqnarray*}
It is easy to check that $\pi|_{vv}> \tilde{\pi}|_{\pi_{it}|_{vv}}, \tilde{\pi}_{i}|_{\pi_i|_{vv}}, \tilde{\pi}_{i}|_{  \pi_{it}|_{vv}}$.

If  $\pi|_{vv}=w_1\cdots w_{t-1}vvw_{s+1}\cdots w_{m}$ and  $v=w_tw_{t+1}\cdots w_{s}, s>1$, then
$$
\llbracket\pi|_{vv}\rrbracket'_{vv}=[[\pi|_{vv}]^{(1)}]'_{[v]^{(1)}[v]^{(1)}}.
$$
Based on  the new letters $[w_1], [w_2],\cdots, [w_{m}]$ and  Theorem \ref{th2.17},   we have
$$
[[\pi|_{vv}]^{(1)}]'_{[v]^{(1)}[v]^{(1)}}=\alpha^{(1)} \pi^{(1)}|_{[[v]^{(1)}[v]^{(1)}]}+\sum \alpha^{(1)}_i \pi^{(1)}_i|_{[[v]^{(1)}[v]^{(1)}]},
$$
Since  $[[v]^{(1)}[v]^{(1)}]=[vv]$ and  $\overline{[w_i]}=w_i$ for $i=1,2,\cdots, m$,    we can get
$$
\pi^{(1)}|_{[[v]^{(1)}[v]^{(1)}]}=\tilde{\alpha}\pi|_{[vv]}+ \sum \alpha^{(j)} \pi_{1j}|_{[vv]}
$$
and
$$
\pi^{(1)}_i|_{[[v]^{(1)}[v]^{(1)}]}=\sum \alpha_i^{(j)} \pi_{ij}|_{[vv]}
$$
with $\pi|_{vv}> \pi_{1j}|_{vv},  \pi_{ij}|_{vv}$. Therefore,
$$
\llbracket\pi|_{vv}\rrbracket'_{vv}=\alpha^{(1)} \tilde{\alpha}\pi|_{[vv]}+ \sum \alpha^{(1)}\alpha^{(j)} \pi_{1j}|_{[vv]}+\sum_{i,j} \alpha^{(1)}_i\alpha_i^{(j)} \pi_{ij}|_{[vv]}.
$$
The proof is complete. \hfill $\square$

Let $u, v\in \textbf{OSLSW}(X)$ and $v$ is a subword of $u$. We will define $\llbracket u\rrbracket_{v}$ as follows. There are four cases to consider.

Case 1. If $u, v$ are  pure operated super-Lyndon-Shirshov words, then we  have defined
$
\llbracket u\rrbracket_{v}
$
by  Theorem \ref{th5.3}.

Case 2.  If $u$ is a   pure operated super-Lyndon-Shirshov word and $v=v_1v_1$ with $|v_1|=1$, then we  define
$\llbracket u\rrbracket_{v}=\llbracket u\rrbracket'_{v_1v_1}.$

Case 3. If $u=u_1u_1$  with $|u_1|=1$ and  $v$ is a  pure operated super-Lyndon-Shirshov word, then by  Lemma \ref{lem5.2}, $u=u_1\pi|_{v}$ or $u=\pi|_{v}u_1$, where   $u_1=\pi|_{v}$.   Define
$
\llbracket u\rrbracket_{v}=[[u_1]\llbracket u_1\rrbracket_{v}]
$
if $u=u_1\pi|_{v}$, and
$
\llbracket u\rrbracket_{v}=[\llbracket u_1\rrbracket_{v}[u_1]]
$
if $u=\pi|_{v}u_1$.

Case 4. If $u=u_1u_1$  and $v=v_1v_1$ with $|u_1|=|v_1|=1$, then by  Lemma \ref{lem5.2}, $u=u_1\pi|_{v}=u_1\pi|_{v_1v_1}$ or $u=\pi|_{v}u_1=\pi|_{v_1v_1}u_1$, where  $u_1=\pi|_{v}=\pi|_{v_1v_1}$.    Define
$
\llbracket u\rrbracket_{v}=[[u_1]\llbracket u_1\rrbracket'_{v_1v_1}]
$
if $u=u_1\pi|_{v}=u_1\pi|_{v_1v_1}$,  and $
\llbracket u\rrbracket_{v}=[\llbracket u_1\rrbracket'_{v_1v_1}[u_1]]
$
if $u=\pi|_{v}u_1=\pi|_{v_1v_1}u_1$.

\begin{theorem}   Let $u, v\in \textbf{OSLSW}(X)$ and $v$ is a subword of $u$. Then $\overline{\llbracket u\rrbracket_{v}}=u$.
\end{theorem}
{\bf Proof.}  By Theorems \ref{th5.3} and \ref{th5.4}, we can get the result. \hfill $\square$
\begin{definition}
Let $\pi\in \langle P;  X\rangle^\star$ and $f\in \textbf{OLS}(X) \subseteq \textbf{OAS}(X)$.   If $\pi|_{\bar{f}} \in \textbf{OSLSW}(X)$, then
$$
\llbracket\pi|_{_{f}}\rrbracket_{_{\bar{f}}}=\llbracket\pi|_{_{\bar{f}}}\rrbracket_{_{\bar{f}}}|_{_{[\bar{f}]\mapsto f}}
$$
is called a special normal $f$-word.
\end{definition}

\begin{proposition}\label{co2.14}
Let $f\in \textbf{OLS}(X)$ be homogenous, $\pi\in \langle P;  X\rangle^\star$  and $\pi|_{\bar{f}}\in \textbf{OSLSW}(X)$. Then we have
$$
\llbracket\pi|_{_{f}}\rrbracket_{_{\bar{f}}}=\alpha_1\pi|_{_{f}}+\sum\alpha_i \pi_i|_{_{f}},
$$
where $\alpha_1=\frac{lc(\llbracket\pi|_{_{ f}}\rrbracket_{_{\bar{f}}})}{lc(f)}$,  each $\alpha_i\in k$, $\pi_i\in \langle P;  X\rangle^\star$  and $\pi_i|_{_{\bar{f}}} <_{_{Dl}}  \pi|_{_{\bar{f}}}$. It follows that $\overline{\llbracket\pi|_{_{f}}\rrbracket_{_{\bar{f}}}}=\pi|_{_{\overline{f}}}$.
\end{proposition}
{\bf Proof.} It is easy to check    that the result is true  by Theorems \ref{th5.3} and \ref{th5.4}. \hfill $\square$

\subsection{Composition-Diamond lemma for  operated  Lie  superalgebras}
In this subsection,  we establish  the  Gr\"{o}bner-Shirshov bases  theory for  operated  Lie  superalgebras.  Especially, we prove    Composition-Diamond lemma for  operated  Lie  superalgebras.

Let $f, g \in \textbf{OLS}(X)$ be  homogenous.  Define  three kinds of compositions as follows.
\begin{enumerate}
\item[(a)]\
If  $w=\bar{f}a=b\bar{g}$ is a pure  operated  super-Lyndon-Shirshov word, where $a,b\in \langle P;  X\rangle$ with  $bre(w)<bre(\bar{f})+bre(\bar{g})$, then
$$
\langle f,g\rangle_w=lc(\llbracket fa\rrbracket_{_{\overline{f}}})^{-1}\llbracket fa\rrbracket_{_{\overline{f}}} - lc(\llbracket bg\rrbracket_{_{\overline{g}}})^{-1} \llbracket bg\rrbracket_{_{\overline{g}}}
$$
is called the intersection composition of $f$ and $g$ with respect to the ambiguity $w$.

\item[(b)]\  If $w=\overline{f}=\pi|_{\overline{g}}$, where $\pi\in \langle P;  X\rangle^\star$,    then
$$
\langle f,g\rangle_w=lc(f)^{-1}f-lc(\llbracket\pi|_{g}\rrbracket_{_{\overline{g}}})^{-1}\llbracket \pi|_{g}\rrbracket_{_{\overline{g}}}
$$
is called the inclusion composition of $f$ and $g$ with respect to  the ambiguity  $w$.

\item[(c)] If $w=\overline{f}=vv$  with $|v|=1$, then
$ \langle f, f\rangle_w=(f[v])$ is called the multiplication   composition of $f$ with respect to  the ambiguity $w$.
\end{enumerate}

Note that for any $S\subseteq \textbf{OLS}(X)$, there is  a homogenous  subset  $\widetilde{S}$ of $\textbf{OLS}(X)$ such that $Id(S)=Id(\widetilde{S})$, where $Id(S)$ is the  operated superideal of $\textbf{OLS}(X)$ generated by $S$.

If $S$ is a homogenous  subset of $\textbf{OLS}(X)$ and  $h\in \textbf{OLS}(X)$, then $h$ is called trivial modulo $(S, \leq_{_{Dl}})$ if
$$
h=\sum \alpha_i \llbracket \pi_i|_{s_i}\rrbracket_{_{\overline{s_i}}},
$$
where each $\alpha_i\in k,
  \ s_i\in S$, $\llbracket\pi_i|_{s_i}\rrbracket_{_{\overline{s_i}}}$  is a special  normal $s_i$-word and $\pi_i|_{_{\overline{s_i}}}\leq_{_{Dl}} \overline{h}$. If
this is the case, then we write
$$
h\equiv 0\  mod(S, \leq_{_{Dl}}).
$$

\begin{definition}
A homogenous set $S$ is called a
Gr\"{o}bner-Shirshov basis  in   $\textbf{OLS}(X)$ if any
composition $\langle f,g\rangle_w$ of  $f, g\in S$ is trivial modulo $(S, \leq_{_{Dl}})$.
\end{definition}

\begin{lemma}\label{le3.6}
Let  $S$ be a homogenous subset of $\textbf{OLS}(X)$  and
$$
Irr(S):= \{[w]|w\in \textbf{OSLSW}(X), w\neq \pi|_{\bar{s}},  s\in S,  \pi\in\langle P;  X\rangle^\star\}.
$$
Then, for any  $h\in \textbf{OLS}(X)$, $h$ can be expressed as
$$
h=\sum\alpha_i[u_i]+
\sum\beta_j\llbracket\pi_j|_{s_j}\rrbracket_{_{\overline{s_j}}},
$$
where each $u_i\in \textbf{OSLSW}(X), u_i\leq_{_{Dl}} \bar{h}$, $s_i\in S$, $\llbracket \pi_j|_{s_j}\rrbracket_{_{\overline{s_j}}}$ is a special  normal $s_j$-word and $\pi_j|_{_{\overline{s_j}}}\leq_{_{Dl}} \bar{h}$. \hfill $ \square$\\
\end{lemma}
{\bf Proof.} Since $h\in \textbf{OLS}(X)$,
$
h=\sum  \alpha_{i}[u_{i}],
$
where each $u_{i}\in \textbf{OSLSW}(X)$ and  $u_{i}>_{_{ Dl}} u_{i+1}$.
If $[u_1]\in Irr(S)$, then let
$
h_{1}=h-\alpha_{1}[u_1].
$
Otherwise,   there exists $s_1\in{S}$ such that $u_1=\pi|_{_{\overline{s_1}}}$. Let
$
h_1=h-\alpha_1lc([u_1])lc( \llbracket\pi|_{s_1}\rrbracket_{_{\overline{s_1}}})^{-1}  \llbracket\pi|_{s_1}\rrbracket_{_{\overline{s_1}}}.
$
In both of  the above two cases, we have  $h_1\in \textbf{OLS}(X)$ and $h_1<_{_{Dl}} \bar{h}$. Then, by induction on $\bar{h}$,  we can obtain the result.
 \hfill $ \square$\\

\begin{proposition}\label{l4.2}\label{lemm3.7}
Let $S$ be a Gr\"{o}bner-Shirshov  basis  in  $\textbf{OLS}(X)$, $\pi_1,\pi_2\in  \langle P; X\rangle^\star$
 and $s_1, s_2\in S$. If
 $w=\pi_1|_{\overline{s_1}}=\pi_2|_{\overline{s_2}}$ and $\alpha_1=lc(s_1)^{-1}, \ \alpha_2=lc(s_2)^{-1}$, then
$$
\pi_1|_{\alpha_1s_1}-\pi_2|_{\alpha_2s_2}=\sum \delta_j \tau_j|_{s_j},
$$
where   $ \delta_j\in k$, $\tau_j\in  \langle P; X\rangle^\star, s_j\in S$ and $\tau_j|_{\overline{s_j}}<_{_{Dl}}w$.
\end{proposition}
{\bf Proof:} \ There are three cases to consider.

(I)\ \ The  words $\overline{s_1}$ and $\overline{s_2}$ are
disjoint in $w$. Then there exits a
$(\star_1,\star_2)$-word $\Pi$ such that
$$\Pi|_{\overline{s_1},\
\overline{s_2}}=\pi_1|_{\overline{s_1}}=\pi_2|_{\overline{s_2}}.
$$
A simple computation gives
\begin{eqnarray*}&& \pi_1|_{\alpha_1 s_1}-\pi_2|_{\alpha_2 s_2}\\
&=& \Pi|_{\alpha_1 s_1, \ \overline{s_2}}-\Pi|_{\overline{s_1}, \  \alpha_2 s_2}\\
&=&(-\Pi|_{\alpha_1 s_1,\ \alpha_2s_2-\overline{s_2}}+\Pi|_{\alpha_1s_1, \ \alpha_2s_2})+(\Pi|_{\alpha_1s_1-\overline{s_1}, \ \alpha_2 s_2}-\Pi|_{\alpha_1s_1, \
\alpha_2s_2})\\
&=&-\Pi|_{\alpha_1 s_1,\ \alpha_2s_2-\overline{s_2}}+\Pi|_{\alpha_1s_1-\overline{s_1}, \ \alpha_2 s_2}
\end{eqnarray*}
It is worth pointing out that $\overline{\alpha_1s_1-\overline{s_1}}<_{_{Dl}}\overline{s_1}$.  Then, we can assume that
$$
\alpha_1s_1-\overline{s_1}=\sum \beta_t u_t,
$$
where $\beta_t\in k$ and  $u_t<_{Dl}\overline{s_1}$.   Thus
$$
\Pi|_{\alpha_1s_1-\overline{s_1}, \ \alpha_2 s_2}=\sum \beta_t\alpha_2 \Pi|_{u_t, \  s_2}=\sum\alpha_{2_t}\pi_{2_t}|_{s_2},
$$
where $\alpha_{2_t}=\beta_t\alpha_2$ and $\pi_{2_t}|_{s_2}=\Pi|_{u_t,     s_2}$.   Since $u_t<_{Dl}\overline{s_1}$ for each $t$,  it follows that $\pi_{2_t}|_{\overline{s_2}}<w$.  Similarly, we have
$$
\Pi|_{\alpha_1 s_1,\ \alpha_2s_2-\overline{s_2}}=\sum\alpha_{1_l}\pi_{1_l}|_{s_1}
$$
where for each $l$, $\pi_{1_l}|_{ \overline{s_1} }<w$.  Therefore, we can get
$$
 \pi_1|_{\alpha_1s_1}- \pi_2|_{\alpha_2s_2}=\sum\alpha_{2_t}\pi_{2_t}|_{s_2}+\sum\alpha_{1_l}\pi_{1_l}|_{s_1},
$$
where each  $\pi_{2_t}|_{\overline{s_2}}, \
\pi_{1_l}|_{\overline{s_1}}<_{_{Dl}}w$.\\

(II) The  words  $\overline{s_1}$ and
 $\overline{s_2}$ have nonempty intersection in $w$ but do not
 include each other, say $\overline{s_1}=e_1e_2$ and $\overline{s_2}=e_2e_3$, where $e_1, e_2, e_3$ are nonempty word.

Case II-1: If $\overline{s_1}\neq \overline{s_2}$  and $\overline{s_1}a=e_1e_2e_3=b\overline{s_2}$,  then  by Proposition \ref{pro5.1},  $\overline{s_1}a=e_1e_2e_3=b\overline{s_2}$   is  a  pure  operated  super-Lyndon-Shirshov word.  It is easy to see that there  exists a
$\star$-word $\Pi$ such that
$$
\Pi|_{\overline{s_1}a}=\pi_1|_{\overline{s_{_1}}}=\pi_2|_{\overline{s_2}}=\Pi|_{b\overline{s_2}}.
$$
By Proposition  \ref{co2.14}, we can get
$$
lc(\llbracket s_{_1}a\rrbracket_{\overline{s_1}})^{-1}\llbracket s_{_1}a\rrbracket_{\overline{s_1}}= \alpha_1s_{_1}a+\sum \beta_{1i} c_{1i}s_{_1}d_{1i},
$$
$$
lc(\llbracket bs_{_2}\rrbracket_{\overline{s_1}})^{-1}\llbracket bs_{_2}\rrbracket_{\overline{s_2}}= \alpha_2bs_{_2}+\sum \beta_{2i} c_{2i}s_{_2}d_{2i},
$$
where each $\beta_{1i}, \beta_{2i}\in k$, $c_{1i}\overline{s_{_1}}d_{1i}<_{_{Dl}}\overline{s_{_1}}a$ and $c_{2i}\overline{s_{_2}}d_{2i}<_{_{Dl}}  b\overline{s_{_2}}$.
It follows that
\begin{eqnarray*}
&&\pi_1|_{\alpha_1 s_{_1}}-\pi_{_2}|_{\alpha_2 s_{_2}}=\Pi|_{\alpha_1 s_{_1}a}-\Pi|_{\alpha_2 bs_{_2}}\\
&=&\Pi|_{lc(\llbracket s_{_1}a\rrbracket_{\overline{s_1}})^{-1}\llbracket s_{_1}a\rrbracket_{\overline{s_1}}-lc(\llbracket bs_{_2}\rrbracket_{\overline{s_1}})^{-1}\llbracket bs_{_2}\rrbracket_{\overline{s_2}}}-\sum \beta_{1i}\Pi|_{ c_{1i}s_{_1}d_{1i}} + \sum \beta_{2i} \Pi|_{ c_{2i}s_{_2}d_{2i}}\\
&=& \Pi|_{ \langle s_1, s_2 \rangle_{\overline{s_1}a}}-\sum \beta_{1i}\Pi|_{ c_{1i}s_{_1}d_{1i}} + \sum \beta_{2i} \Pi|_{ c_{2i}s_{_2}d_{2i}}.
\end{eqnarray*}
Since  $S $ is a Gr\"{o}bner-Shirshov basis in  $\textbf{OLS}(X)$, by Propositon \ref{co2.14},
 we have
$$
\langle s_1, s_2 \rangle_{\overline{s_1}a}= \sum\tilde{\alpha}_j\llbracket\tilde{\tau}_j|_{\tilde{s}_j}\rrbracket_{\overline{\tilde{s}_j}}=\sum \gamma_n\varphi_n|_{s_n},
$$
where each $\gamma_n\in k, \ \varphi_n\in \langle P; X\rangle^\star, \ s_n\in
S$ and $\varphi_n|_{\overline{s_n}}\leq_{_{Dl}}\overline{\langle s_1, s_2 \rangle_{\overline{s_1}a}} <_{_{Dl}}\overline{s_1}a$. Let $
\Pi|_{\varphi_n|_{s_n}}=\Pi_{n}|_{s_n}, \Pi|_{ c_{1i}s_{_1}d_{1i}}=\Pi_{1i}|_{s_1}, \Pi|_{ c_{2i}s_{_2}d_{2i}}=\Pi_{2i}|_{s_2}.$
Then
$$
\alpha_1\pi_1|_{ s_{_1}}-\alpha_2\pi_{_2}|_{ s_{_2}}=\sum \gamma_n\Pi_{n}|_{s_n}+\sum \beta_{1i}\Pi_{1i}|_{s_1}+\sum \beta_{2i}\Pi_{2i}|_{s_2}
$$
where each $\Pi_{n}|_{\overline{s_n}}, \Pi_{1i}|_{\overline{s_1}}, \Pi_{2i}|_{\overline{s_2}} < w$.\\

Case II-2:  If $\overline{s_1}=\overline{s_2}$, then    $\overline{s_1}=\overline{s_2}=uu$ with $|u|=1$.  Otherwise,  $\overline{s_1}=\overline{s_2}$ is a pure  operated  super-Lyndon-Shirshov word. It follows that
$$
\overline{s_1}=e_1e_2=e_2e_3\succ_{lex}e_2e_1,
$$
which gives that $e_3\succ_{lex}e_1$.  Thus, $e_3e_2\succ_{lex}e_1e_2=e_2e_3=\overline{s_1}$, a contradiction.  By Proposition \ref{pro5.1}, we have $e_1=e_2=e_3=u$.    Moreover,  we have
\begin{eqnarray*}
\pi_1|_{\alpha_1 s_1}-\pi_2|_{\alpha_2 s_2}
= \Pi|_{\alpha_1 s_1u -\alpha_2us_2}.
\end{eqnarray*}
By Lemma \ref{lem2.7}, we have $[u]=u+ H(u)$, where $H(u)=\sum \gamma_tu_t$ with $u>_{Dl}u_t$ and $bre(u)=bre(u_t).$ For $i=1,2$,  let $\alpha_i s_i=\overline{s_i}+H(\alpha_i s_i)$, where $H(\alpha_i s_i)=\sum \beta_{ij}w_{ij}$ with $\overline{s_i}>_{Dl}w_{ij}$,
Consider  the fact that
\begin{eqnarray*}
 &&\langle \alpha_1 s_1, \alpha_1 s_1\rangle_{uu}\\
 &=&(\alpha_1 s_1[u])=\alpha_1 s_1[u]-[u]\alpha_1 s_1\\
&=& \alpha_1 s_1(u+ H(u))-(u+ H(u))\alpha_1 s_1\\
&=&  \alpha_1 s_1u-u\alpha_1 s_1+  \alpha_1 s_1H(u)-  H(u) \alpha_1 s_1\\
&=&  \alpha_1 s_1u-u(\overline{s_1}+H(\alpha_1 s_1))+  \alpha_1 s_1H(u)-  H(u) \alpha_1 s_1\\
&=&  \alpha_1 s_1u-u\overline{s_2}-uH(\alpha_1 s_1) +  \alpha_1 s_1H(u)-  H(u) \alpha_1 s_1\\
&=&  \alpha_1 s_1u-u\alpha_2 s_2+ uH(\alpha_2 s_2)- uH(\alpha_1 s_1) +  \alpha_1 s_1H(u)-  H(u) \alpha_1 s_1.
\end{eqnarray*}
By an easy computation, we can get
\begin{eqnarray*}
H&=&uH(\alpha_2 s_2)- uH(\alpha_1 s_1) +  \alpha_1 s_1H(u)-  H(u) \alpha_1 s_1\\
&=&  \sum \beta_{2j}uw_{2j}-\sum \beta_{1j}uw_{1j}+\sum \gamma_t\overline{s_1}u_t\\
&&+\sum_{t,j}\gamma_t\beta_{1j}w_{1j}u_t+\sum \gamma_t u_t\overline{s_1}+\sum_{t,j}\gamma_t\beta_{1j}u_tw_{1j}.\\
\end{eqnarray*}
Since  $u>_{Dl}u_t$, $\overline{s_i}>_{Dl}w_{ij}, i=1,2$,  it follows that
$$
uw_{2j}, uw_{1j},  \overline{s_1}u_t, w_{1j}u_t,  u_t\overline{s_1},  u_tw_{1j} \leq_{Dl} \overline{\langle \alpha_1 s_1, \alpha_1 s_1\rangle_{uu}}<_{Dl} uuu.
$$ To   simplify  notation, we write  $H=\sum \delta_lv_l,$  where  $v_l\leq_{Dl}\overline{\langle \alpha_1 s_1, \alpha_1 s_1\rangle_{uu}}<_{Dl}uuu$.   It follows that
$$
\alpha_1 s_1u-u\alpha_2 s_2  =\langle \alpha_1 s_1, \alpha_1 s_1\rangle_{uu}-\sum\delta_lv_l.
$$
which gives that
$$
\pi_1|_{\alpha_1 s_1}-\pi_2|_{\alpha_2 s_2}=\Pi|_{\alpha_1 s_1u -\alpha_2us_2}=\Pi|_{\langle \alpha_1 s_1, \alpha_1 s_1\rangle_{uu}}- \sum \delta_l \Pi|_{v_l}.
$$
Since  $S $ is a Gr\"{o}bner-Shirshov basis in  $\textbf{OLS}(X)$,  similar to the   Case II-1, we can get the result.

(III) One of  words $\overline{s_1}$,
 $\overline{s_2}$ is contained in the other one.
For example, let $ \overline{s_1}=\pi|_{\overline{s_2}}$ for some
$\star$-word $\pi$. Then we have
$$
w=\pi_2|_{\overline{s_2}}=\pi_1|_{\pi|_{\overline{s_2}}}.
$$
It follows that
$$
\pi_1|_{ \alpha_1 s_1}-\pi_2|_{ \alpha_2 s_2}=\pi_1|_{ \alpha_1 s_1}-\pi_1|_{\pi|_{\alpha_2 s_2}}=\pi_1|_{
\alpha_1s_1-\alpha_2\pi|_{s_2}}.
$$
Since  $S $ is a Gr\"{o}bner-Shirshov basis in  $\textbf{OLS}(X)$, similar to the proof  of  Case II-1, we  can obtain the result. \hfill
$\square$
\\

\begin{theorem}\label{cdll}
(Composition-Diamond lemma for operated  Lie  superalgebras)\ Let $S$ be a  homogenous subset  of $\textbf{OLS}(X)$  and $Id(S)$ the operated superideal of $\textbf{OLS}(X)$ generated by $S$.
 Then the following statements are equivalent:
\begin{enumerate}
\item[(a)]The set  $S $ is a Gr\"{o}bner-Shirshov basis  in $\textbf{OLS}(X)$.
\item[(b)] If  $f\in Id(S)$, then $\bar{f}=\pi|_{\bar{s}}$ for some $s\in S$ and $\pi\in\langle P;  X\rangle^\star$.
\item[(c)] The set
$$
Irr(S)= \{[w]|w\in \textbf{OSLSW}(X), w\neq \pi|_{\bar{s}}, s\in S,  \pi\in\langle P;  X\rangle^\star \}
$$
is a  linear basis of   $\textbf{OLS}(X|S)=\textbf{OLS}(X)/Id(S)$.
\end{enumerate}
\end{theorem}
{\bf Proof.}
 (a)$\Longrightarrow$ (b)\ \ If   $0\neq f\in
Id(S)$, then it is easy to see that  $f\in \textbf{Id}(S)$, where $\textbf{Id}(S)$ is the  associative operated  superideal of the free associative operated superalgebra $\textbf{OAS}( X)$   generated by $S$.  It follows that
$$
f= \sum\limits_{i=1}^{n} \alpha_i\pi_i|_{s_i},
$$
where each $\alpha_i\in k$,  $\pi_i\in \langle P; X\rangle^\star$ and  $s_i\in
S$.

Let $w_i= \pi_i|_{_{\overline{ s_i}}}$. We arrange these leading words in nonincreasing order by
$$
w_1= w_2=\cdots=w_m >_{_{Dl}}w_{m+1}\geq_{_{Dl}} \cdots\geq_{_{Dl}} w_n.
$$
We prove the result by  induction  on $w_1$ and $m$. It is worth  pointing out  that $w_1\geq \overline{f}$ and $m\geq 1$.

If  $w_1=\overline{ f}$ or  $m=1$, then it is  easy to see that the result is true.

Assume that the result is true  for any $(\widetilde{w}_1, \widetilde{m})<(w_1, m)$.    We can also suppose that   $w_1>\bar{f}$ and $m>1$.  Note that
$\pi_1|_{\overline{s_1}}=w_1=w_2=\pi_2|_{\overline{s_2}}$.
 As  $S $ is a Gr\"{o}bner-Shirshov
basis in $\textbf{OLS}(X)$ we have
$$
lc(s_1)^{-1}\pi_1|_{s_1}-lc(s_2)^{-1}\pi_2|_{s_2}= \sum \delta_i  \tau_j|_{s_j}
$$
by Proposition \ref{lemm3.7},  where $ \delta_i\in k$ and $\tau_j|_{\overline{s_j}}<_{_{Dl}}w_1$. An easy computation gives that
\begin{eqnarray*}
&& \alpha_1\pi_1|_{s_1}+ \alpha_2\pi_2|_{s_2}\\
&=&(\alpha_1+\alpha_2lc(s_2) lc(s_1)^{-1} )\pi_1|_{s_1}
 +\alpha_2lc(s_2) (lc(s_2)^{-1} \pi_2|_{s_2} - lc(s_1)^{-1} \pi_1|_{s_1} ).
\end{eqnarray*}
It follows that
$$
f=(\alpha_1+\alpha_2lc(s_2) lc(s_1)^{-1})\pi_1|_{s_1}
 +
\sum  \alpha_2lc(s_2)\delta_i  \tau_j|_{s_j}
+
\sum\limits_{i=3}^{n} \alpha_i\pi_i|_{s_i}.
$$
There are two cases to consider: (i) $m=2$ and $\alpha_1+\alpha_2lc(s_2) lc(s_1)^{-1}=0$,   and (ii) either $m>2$ or $\alpha_1+\alpha_2lc(s_2) lc(s_1)^{-1}\neq 0$. On both cases, we can obtain the result by induction.

 $(b)\Rightarrow(c)$   \ Suppose that
$
\sum\alpha_i[u_i]=0
$
in $\textbf{OLS}(X|S)$, where each $[u_i]\in Irr(S)$ and  $u_i>_{_{  Dl}}  u_{i+1}$.  That
is,
$
\sum\alpha_i[u_i]\in{\textbf{Id}(S)}.
$
Then each $\alpha_i$ must be 0. Otherwise, say $\alpha_1\neq0$,
since
$
\overline{\sum\alpha_i[u_i]}=u_1
$
and by (b), we have  $u_1\notin Irr(S)$,    a contradiction. Therefore, $Irr(S)$ is linear independent. By Lemma \ref{le3.6},   $Irr(S)$
is a  linear basis of     $\textbf{OLS}(X|S)$.

$(c)\Rightarrow(a)$  \ For any composition $\langle
f,g\rangle_w$  of $f,g\in{S}$, we have $\langle
f,g\rangle_w\in{\textbf{Id}(S)}$. Then, by (c) and   Lemma \ref{le3.6},
$
 \langle
f,g\rangle_w=\sum\beta_j\llbracket \pi_j|_{s_j}\rrbracket_{_{\overline{s_j}}},
$
where each $\beta_j\in k,   \pi_j\in\langle P;  X\rangle^\star, s_j\in S, \pi_j|_{_{\overline{s_j}}} \leq_{_{ Dl}} \overline{\langle f,g\rangle_w}$. This proves that $S$ is a
Gr\"{o}bner-Shirshov basis in $\textbf{OLS}(X)$.   \hfill $ \square$\\

\begin{theorem}
Let $S\subseteq \textbf{OLS}(X)\subseteq \textbf{OAS}(X)$ be    homogenous.   Then the following two statements are equivalent.
\begin{enumerate}
\item[(a)] The set  $S $ is a Gr\"{o}bner-Shirshov basis  in $\textbf{OLS}(X)$.
\item[(b)] The set  $S $ is a Gr\"{o}bner-Shirshov basis  in $\textbf{OAS}(X)$.
\end{enumerate}
\end{theorem}
{\bf Proof.} By Theorems \ref{cdla} and \ref{cdll}, we can get the result easily.   \hfill $ \square$\\

\section{Construction of  free Lie Rota-Baxter  superalgebras}

\subsection{A  Gr\"{o}bner-Shirshov basis of  a free  Lie Rota-Baxter  superalgebra}
In this  section, as an application of Theorem \ref{cdll}, we find a  Gr\"{o}bner-Shirshov basis of  a  free Lie  Rota-Baxter superalgebra. As consequently,  we can obtain      a linear basis of  a free Lie  Rota-Baxter superalgebra  by  Composition-Diamond lemma for operated  Lie superalgebras.\\

\begin{definition}(\cite{amm2015,  whc2010})
A Lie superalgebra $L=L_{\bar{0}}\oplus L_{\bar{1}}$  equipped with an  even linear map $ P:L\rightarrow L$ satisfying
$$
[P(x)P(y)]  = P([P(x) y] )-(-1)^{|x||y|}P([P(y)x])+\lambda P ( [x y ])
$$
for all homogenous elements $x$ and $y$,   is called a  Lie   Rota-Baxter  superalgebra of weight $\lambda$ or Lie $\lambda$-Rota-Baxter  superalgebra.
\end{definition}

It is obvious  that any  Lie $\lambda$-Rota-Baxter  superalgebra   is an operated  Lie  superalgebra.

Let   $\textbf{OLS}(X)$   be the free  operated  Lie    superalgebra   on  the  $\mathbb{Z}_2$-graded set $X=X_{\bar{0}}\cup X_{\bar{1}}$.   It is easily  seen that
$$
 \textbf{RBLS}(X)=\textbf{OLS}(X|\widetilde{S})=\textbf{OLS}(X)/\textbf{Id}(\widetilde{S})
$$
is a free Lie $\lambda$-Rota-Baxter superalgebra on  the $\mathbb{Z}_2$-graded set $X=X_{\bar{0}}\cup X_{\bar{1}}$, where $\textbf{Id}(\widetilde{S})$ is a operated Lie super-ideal  of $\textbf{OLS}(X)$ generated by
$$
\widetilde{S}=  \left\{
\left. \begin{array}{ll}
[P([u])P([v])]-P([P([u])  [v]] ) \\
\ \ \   \ \ \ +(-1)^{|u||v|}P([P([v])[u]])-\lambda P ( [[u]  [v]]))
\end{array}\right| u, v\in \textbf{OSLSW}(X)
 \right\}.
$$
Let us  denote by $S$    the  set consisting of the following polynomials in $\textbf{OLS}(X)$:
$$
(P([w])P([w]))- 2P((P([w])[w]))  -\lambda P(([w][w])),  \ |w|=1,
$$
$$
(P([u])P([v]))- P((P([u])[v]))+(-1)^{|u||v|}P((P([v])[u]))-\lambda P(([u][v])),
$$
where $u, v,  w \in  \textbf{OSLSW}(X)$ and  $u>_{_{Dl}}v$.
It is easy to check that  $\textbf{Id}(S)=\textbf{Id}(\widetilde{S})$. It follows that
$$
\textbf{RBLS}(X)=\textbf{OLS}(X|S)=\textbf{OLS}(X)/\textbf{Id}(S).
$$

\begin{theorem}\label{th4.1}With the order  $>_{_{Dl}}$  on $\langle P; X\rangle$,  the set $S$ is a Gr\"{o}bner-Shirshov basis in $\textbf{OLS}(X)$. It follows that the set
$$
Irr(S)=  \left\{
[z]\in  \textbf{OSLSM}(X)
 \left|
 \begin{array}{ll}
 z\neq \pi|_{P(u)P(v)}, \ \  u >_{_{Dl}} v, \\
   z\neq \pi|_{P(w)P(w)},  \   |w|=\bar{1},\\
 u, v, w\in \textbf{OSLSW}(X), \pi\in  \langle P; X\rangle^\star
\end{array}
\right. \right\}
$$
is a linear basis of the free  Lie $\lambda$-Rota-Baxter  superalgebra $\textbf{RBLS}(X)$ on  the  $\mathbb{Z}_2$-graded set $X=X_{\bar{0}}\cup X_{\bar{1}}$.
\end{theorem}

\subsection{Proof of Theorem \ref{th4.1} }
In this subsection, we will prove  Theorem \ref{th4.1}. Let us  denote by
$$
g_{w}=(P([w])P([w]))- 2P((P([w])[w]))  -\lambda P(([w][w])),  \ |w|=1,
$$
$$
f_{u, v} =(P([u])P([v]))- P((P([u])[v]))+(-1)^{|u||v|}P((P([v])[u]))-\lambda P(([u][v])),
$$
where $u, v,  w \in \textbf{OSLSW}(X)$ with $u>_{_{Dl}}v$.\\

(1)\   We will check that all the multiplication compositions are trivial. That is to check that   for any $w\in \textbf{OSLSW}(X)$ with $|w|=1$,  $$\langle g_w, g_w\rangle_{P(w)P(w)} \equiv 0\  {mod}(S, \leq_{_{Dl}}).$$  Note that $|P([w])|=|[w]|=|w|=1$. By the  super-Jacobi identity and super-skew-symmetry,  we have
 \begin{eqnarray*}
(P([w])((P([w])[w]))&=& ((P([w])P([w]))[w])+(-1)^{|P([w])||P([w])|} (P([w](P([w])[w])))\\
&=&  ((P([w])P([w]))[w])-(P([w])((P([w])[w])),\\
((P([w])[w])[w])&=&(P([w])([w][w]))-(-1)^{|P([w])||[w]|}([w](P([w])[w])) \\
&=& (P([w])([w][w]))-((P([w])[w])[w]),
 \end{eqnarray*}
which gives that
$$
(P([w])((P([w])[w]))
=\frac{1}{2}  ((P([w])P([w]))[w]), \  (([w][w])P([w]))=-2((P([w])[w])[w]).
$$
It follows that
 \begin{eqnarray*}
&&\langle g_w, g_w\rangle_{P(w)P(w)}\\
&=& (((P([w])P([w]))- 2P((P([w])[w]))  -\lambda P(([w][w])))P([w]))\\
&=& - 2(P((P([w])[w]))P([w]))  -\lambda (P(([w][w]))P([w]))\\
&\equiv& -2P((P((P([w])[w]))[w])+((P([w])[w])P([w]))+\lambda((P([w])[w])[w]))\\
&&-\lambda P((P(([w][w]))[w])+ (([w][w])P([w]))+\lambda (([w][w])[w]))\\
&\equiv&P((((P([w])P([w])) -2P((P([w])[w]))    -\lambda P(([w][w]))) [w]))\\
&\equiv& 0\  {mod}(S, \leq_{_{Dl}}).
\end{eqnarray*}

(2)\  All the  possible intersection and inclusion  compositions of
polynomials in $S$ are listed as below:
$$\langle f_{u,v}, f_{v, w}\rangle_{w_1},\ w_1=P(u)P(v)P(w),\ u >_{_{Dl}}v >_{_{Dl}} w,$$
$$\langle f_{\pi|_{P(u)P(v)},w}, f_{u, v}\rangle_{w_2},\ w_2=P(\pi|_{P(u)P(v)})P(w), \  u>_{_{Dl}} v, \ \pi|_{P(u)P(v)} >_{_{Dl}} w,$$
$$\langle f_{u, \pi|_{P(v)P(w)}}, f_{v, w}\rangle_{w_3},\ w_3=P(u)P(\pi|_{P(v)P(w)}),\  v >_{_{Dl}}  w, \ u>_{_{Dl}}\pi|_{P(v)P(w)},$$
$$\langle f_{u,v}, g_{v}\rangle_{w_4},\ w_4=P(u)P(v)P(v),\ u >_{_{Dl}}v, \ |v|=\bar{1},$$
$$\langle f_{\pi|_{P(u)P(u)},v}, g_{u}\rangle_{w_5}, \ w_5=P(\pi|_{P(u)P(u)})P(v), \   \pi|_{P(u)P(u)} >_{_{Dl}}v, |u|=\bar{1},$$
$$\langle f_{u, \pi|_{P(v)P(v)}}, g_{v}\rangle_{w_6},\ w_6=P(u)P(\pi|_{P(v)P(v)}),\    u>_{_{Dl}}\pi|_{P(v)P(v)}, |v|=\bar{1},$$
$$\langle g_{\pi|_{P(u)P(u)}}, g_{u}\rangle_{w_7},\ w_7=P(\pi|_{\overline{g(u)}})P(\pi|_{P(u)P(u)}), \  |\pi|_{P(u)P(u)}|=\bar{1}, \ |u|=\bar{1},$$
$$\langle g_{\pi|_{P(u)P(u)}}, g_{u}\rangle_{w_8},\ w_8=P(\pi|_{P(u)P(u)})P(\pi|_{\overline{g(u)}}), \  |\pi|_{P(u)P(u)}|=\bar{1}, \ |u|=\bar{1},$$
$$\langle g_{\pi|_{P(u)P(v)}}, f_{u,v}\rangle_{w_9},\ w_9=P(\pi|_{\overline{f_{u,v}}})P(\pi|_{P(u)P(v)}),\  u>_{_{Dl}}v,  \ |\pi|_{P(u)P(v)}|=\bar{1},$$
$$\langle g_{\pi|_{P(u)P(v)}}, f_{u,v}\rangle_{w_{10}},\ w_{10}=P(\pi|_{P(u)P(v)})P(\pi|_{\overline{f_{u,v}}}),\  u>_{_{Dl}}v,  \ |\pi|_{P(u)P(v)}|=\bar{1},$$
$$\langle g_{u}, f_{u,v}\rangle_{w_{11}},\ w_{11}=P(u)P(u)P(v),\ u >_{_{Dl}}v, \ |u|=\bar{1}.$$
We have checked  that all the inclusion  compositions above  are trivial. Here, we  just check  $\langle f_{u,v}, g_{v}\rangle_{w_4}\equiv 0\ mod(S, \leq_{_{Dl}})$.   An easy computation   gives that
\begin{eqnarray*}
&&(P((P([u])[v]))P([v]))\\
&=&P((P((P([u])[v]))[v])+ ((P([u])[v])P([v])) +\lambda ((P([u])[v])[v]))\\
&=& P((P((P([u])[v]))[v])+ ((P([u])([v] P([v])) -(P([u]P([v]))[v])+\lambda  ((P([u])[v])[v])),\\
&&(P(([u]P([v])))P([v]))\\
&=& P((P(([u]P([v])))[v]) +(([u]P([v]))P([v]))+\lambda ([u]P([v]))[v]))\\
&=&P((P(([u]P([v])))[v]) +\frac{1}{2}(([u](P([v])P([v]))+\lambda ([u]P([v]))[v])),\\
&&(P(([u] [v] ))P([v]))\\
&=&  P((P(([u] [v] ))[v]) +(([u] [v] )P([v]))+\lambda P((([u][v])[v])),\\
&&(P([u]) P(( P([v])[v]))))\\
&=&  P((P([u])( P([v])[v]))+[u]  (P( P([v])[v])))+\lambda ([u](P([v])[v]))),\\
&&(P([u]) P([v][v]))\\
&=&  P((P([u])([v][v]))+([u]  ( P( [v][v]))) +\lambda ([u]( ([v] [v])))).
\end{eqnarray*}
It follows that
\begin{eqnarray*}
&&\langle f_{u,v}, g_{v}\rangle_{w_4}\\
&=&lc([u])^{-1} lc([v])^{-2}\{(f_{u,v}P([v]))-\frac{1}{2}(P([u])g_{v})\}\\
&\equiv& (f_{u,v}P([v]))-\frac{1}{2}(P([u])g_{v}))\\
&\equiv&  P(((P([u])[v]))P([v]))+P(([u]P([v])))P([v])+\lambda (P(([u] [v] ))P([v]))\\
&&- (P([u])P(( P([v])[v]))-\frac{1}{2}\lambda(P([u]) P(([v][v])))\\
&\equiv& 0\ mod(S, \leq_{_{Dl}}).
\end{eqnarray*}

By the Cases (1) and   (2), we have that    $S$ is a Gr\"{o}bner-Shirshov basis in $\textbf{OLS}(X)$.
By Composition-Diamond lemma for oprated Lie  superalgebras,   $Irr(S)$
 is a linear basis of the free  operated  Lie Rota-Baxter  superalgebra $\textbf{RBLS}(X)$.
\hfill $ \square$\\

\noindent{\bf Acknowledgement}: This work was supported   by  the NSF of Lingnan Normal University (No. ZL1811),   the NNSF of China (No. 12071156, 12171222).

\end{document}